\documentclass[authoryear, 11pt]{imsart}
\RequirePackage[OT1]{fontenc}
\RequirePackage{amsthm,amsmath,natbib}
\RequirePackage{hypernat}
\startlocaldefs
\numberwithin{equation}{section}
\theoremstyle{plain}

\endlocaldefs
\usepackage{psfrag, epsfig, graphicx}
\usepackage{amsfonts,amsmath,latexsym,amssymb}
\usepackage[active]{srcltx}
\usepackage{fullpage}
\usepackage[dvips]{color}

\newtheorem{lemma}{Lemma}
\newtheorem{theorem}{Theorem}
\newtheorem{proposition}{Proposition}
\newtheorem{conjecture}{Conjecture}

\newtheorem{assumption}{Assumption}
\newtheorem{remark}{Remark}
\newtheorem{definition}{Definition}

\newtheorem{assumption*}{Assumption}

\renewcommand{\kappa}{\varkappa}

\newcommand{\rd}{{\rm d}}
\newcommand{\e}{\varepsilon}

\newcommand{\cH}{{\cal H}}
\newcommand{\cI}{{\cal I}}
\newcommand{\cJ}{{\cal J}}
\newcommand{\cK}{{\cal K}}

\newcommand{\cP}{{\cal P}}

\newcommand{\cR}{{\cal R}}

\newcommand{\Bh}{\boldsymbol{h}}

\newcommand{\BL}{\boldsymbol{L}}

\newcommand{\Bv}{\boldsymbol{v}}

\newcommand{\blb}{\boldsymbol{\beta}}
\newcommand{\bld}{\boldsymbol{\delta}}

\newcommand{\bll}{\boldsymbol{\lambda}}
\newcommand{\blg}{\boldsymbol{\mu}}

\newcommand{\bleta}{\boldsymbol{\eta}}
\newcommand{\blL}{\boldsymbol{\Lambda}}
\newcommand{\bB}{\mathbb B}
\newcommand{\bC}{\mathbb C}

\newcommand{\bE}{\mathbb E}
\newcommand{\bF}{\mathbb F}

\newcommand{\bH}{\mathbb H}

\newcommand{\bL}{{\mathbb L}}

\newcommand{\bN}{{\mathbb N}}

\newcommand{\bP}{{\mathbb P}}

\newcommand{\bR}{{\mathbb R}}

\newcommand{\bZ}{{\mathbb Z}}
\newcommand{\mD}{\mathfrak{D}}

\newcommand{\mP}{\mathfrak{P}}

\newcommand{\mh}{\mathfrak{h}}

\newcommand{\mJ}{\mathfrak{J}}

\newcommand{\mH}{\mathfrak{H}}

\newcommand{\ma}{\mathfrak{a}}
\newcommand{\mb}{\mathfrak{b}}

\newcommand{\mz}{\mathfrak{z}}
\newcommand{\ttt}{t}

\newcommand{\epr}{\hfill\hbox{\hskip 4pt
                \vrule width 5pt height 6pt depth 1.5pt}\vspace{0.5cm}\par}

\newcommand{\blh}{\boldsymbol{h}}

\begin{document}
\begin{frontmatter}
\title{Estimation in the convolution structure density model. Part II: adaptation over the scale of anisotropic classes.}
\runtitle{Upper bounds}
\begin{aug}
\author[t1]
{\fnms{O.V.} \snm{Lepski}
\ead[label=e1]{oleg.lepski@univ-amu.fr}}
\author[t1]{\fnms{T.} \snm{Willer}
\ead[label=e2]{thomas.willer@univ-amu.fr}}
\thankstext{t1}{This work has been carried out in the framework of the Labex Archim\`ede (ANR-11-LABX-0033) and of the A*MIDEX project (ANR-11-IDEX-0001-02), funded by the "Investissements d'Avenir" French Government program managed by the French National Research Agency (ANR).}
\runauthor{O.V. Lepski and T. Willer}

\affiliation{Aix Marseille Univ, CNRS, Centrale Marseille, I2M, Marseille, France}

\address{Institut de Math\'ematique de Marseille\\
Aix-Marseille  Universit\'e   \\
 39, rue F. Joliot-Curie \\
13453 Marseille, France\\
\printead{e1}\\
\printead{e2}}
\end{aug}
\begin{abstract}

This paper continues the research started in \cite{LW16}.
In  the framework of the
convolution structure density model on $\bR^d$, we address
the problem of adaptive minimax  estimation
with $\bL_p$--loss over the scale of anisotropic Nikol'skii classes.  We fully characterize the behavior of the minimax risk for different
relationships between regularity parameters and norm indexes in the definitions of the functional class and of the risk. In particular, we show that the boundedness
of the function to be estimated leads to an essential improvement of the asymptotic of the minimax risk.
We prove that the selection
rule proposed in Part I leads to the construction of an optimally or nearly optimally (up to logarithmic factor) adaptive estimator.

\end{abstract}


\end{frontmatter}

\section{Introduction}

In the present paper we will be interested in the adaptive estimation in the convolution structure density model. Our considerations here continue the research started
in \cite{LW16}.

Thus,  we observe i.i.d. vectors $Z_i\in\bR^d, i=1,\ldots, n,$ with a common probability density $\mathfrak{p}$ satisfying the following structural assumption
\begin{equation}
\label{eq:convolution structure model}
\mathfrak{p}=(1-\alpha) f +\alpha [f\star g],\quad f\in\bF_g(R),\;\;\alpha\in[0,1],
\end{equation}
where $\alpha\in[0,1]$ and  $g:\bR^d\to\bR$ are supposed to be known and $f:\bR^d\to\bR$ is the function to be estimated.
Recall that for two functions $f,g\in\bL_1\big(\bR^d\big)$
$$
\big[f\star g\big](x)=\int_{\bR^d}f(x-z)g(z)\nu_d(\rd z),\;\;x\in\bR^d,
$$
and for any $\alpha\in[0,1]$, $g\in\bL_1\big(\bR^d\big)$ and $R> 1$,
$$
\bF_g (R)=\Big\{f\in\bB_{1,d}(R): \;(1-\alpha) f +\alpha [f\star g]\in\mP\big(\bR^d\big)\Big\}.
$$
 Furthermore $\mP\big(\bR^d\big)$ denotes the set of probability densities on $\bR^d$,    $\bB_{s,d}(R)$ is the ball of radius $R>0$ in $\bL_s\big(\bR^d\big):=\bL_s\big(\bR^d,\nu_d\big), 1\leq s\leq\infty$ and  $\nu_d$ is the Lebesgue measure on $\bR^d$. At last, for any $U\in\bL_1\big(\bR^d\big)$ let $\check{U}(\ttt):=\int_{\bR^d}U(x)e^{-i\sum_{j=1}^d x_j t_j}\nu_d(\rd x), \ttt\in\bR^d,$ be the Fourier transform of $U$.

 \vskip0.1cm

The convolution structure density model (\ref{eq:convolution structure model}) will be studied for an arbitrary $g\in\bL_1\big(\bR^d\big)$ and $f\in\bF_g(R)$.
Then, except in the case $\alpha=0$, the function $f$ is not necessarily a probability density.

We want to estimate $f$ using the observations $Z^{(n)}=(Z_1,\ldots,Z_n)$. By estimator, we mean any $Z^{(n)}$-measurable map $\hat{f}:\bR^n\to \bL_p\big(\bR^d\big)$. The accuracy of an estimator $\hat{f}$
is measured by the $\bL_p$--risk
$$
 \cR^{(p)}_n[\hat{f}, f]:=\Big(\bE_f \|\hat{f}-f\|_p^p\Big)^{1/p},\;p\in [1,\infty),
$$
where $\bE_f$ denotes the expectation with respect to the probability measure
$\bP_f$ of the observations $Z^{(n)}=(Z_1,\ldots,Z_n)$.
Also, $\|\cdot\|_p$, $p\in [1,\infty)$, is the $\bL_p$-norm on $\bR^d$. The objective is to
construct an estimator of $f$ with a small $\bL_p$--risk.

\subsection{Adaptive estimation}

Let $\bF$ be  a given subset of $\bL_p\big(\bR^d\big)$. For any estimator $\tilde{f}_n$ define its {\it maximal risk}  by
$
\cR^{(p)}_n\big[\tilde{f}_n; \bF\big]=\sup_{f\in\bF}\cR^{(p)}_n\big[\tilde{f}_n; f\big]
$
and its {\it minimax risk} on $\bF$ is given by
\begin{equation}
\label{eq:minmax-risk}
\phi_n(\bF):=\inf_{\tilde{f}_n}\cR^{(p)}_n\big[\tilde{f}_n; \bF\big].
\end{equation}
Here, the infimum is taken over all possible estimators. An estimator whose maximal risk is bounded by $\phi_n(\bF)$ up to some constant factor is called minimax on $\bF$.

Let $\big\{\bF_\vartheta,\vartheta\in\Theta\big\}$ be a collection of subsets of $\bL_p\big(\bR^d,\nu_d\big)$, where $\vartheta$ is a nuisance parameter which may have a very complicated structure.

The problem of adaptive estimation can be formulated as follows:
{\it is it possible to construct a single estimator $\hat{f}_n$
 which would be  simultaneously minimax on each class
 $\bF_\vartheta,\;\vartheta\in\Theta$, i.e.}
$$
 \limsup_{n\to \infty}\phi^{-1}_n(\bF_\vartheta)\cR^{(p)}_n\big[\hat{f}_n; \bF_\vartheta\big]<\infty,\;\;\forall \vartheta\in\Theta?
$$
We refer to this question as {\it the  problem of minimax adaptive
estimation over  the scale  }
$\{\bF_\vartheta,\;\vartheta\in\Theta \}$.
If such an estimator exists, we will call it optimally adaptive. Using the modern statistical language we call the estimator $\hat{f}_n$
{\it nearly optimally adaptive} if
$$
\limsup_{n\to \infty}\phi^{-1}_\frac{n}{\ln n}(\bF_\vartheta) \cR^{(p)}_n\big[\hat{f}_n; \bF_\vartheta\big]<\infty,\;\;\forall \vartheta\in\Theta.
$$
We will be interested in adaptive estimation over the scale
$$
\bF_\vartheta=\bN_{\vec{r},d}\big(\vec{\beta},\vec{L}\big)\cap\bF_{g,\mathbf{\infty}} (R,Q),\;\;\vartheta=\big(\vec{\beta},\vec{r},\vec{L},R,Q\big),
$$
where $\bN_{\vec{r},d}\big(\vec{\beta},\vec{L}\big)$ is the anisotropic Nikolskii
class, see Definition \ref{def:nikolskii}  below. As it was explained in Part I, the adaptive  estimation over the scale $\big\{\bN_{\vec{r},d}\big(\vec{\beta},\vec{L}\big),\;\big(\vec{\beta},\vec{r},\vec{L}\big)\in(0,\infty)^d\times[1,\infty]^d\times(0,\infty)^d\big\}$
can be viewed as the adaptation to anisotropy  and inhomogeneity of the function to be estimated.
Recall also that
$$
\bF_{g,\mathbf{\infty}} (R,Q):=\big\{f\in\bF_g(R): \;(1-\alpha) f +\alpha [f\star g]\in\bB_{\mathbf{\infty},d}(Q)\big\},
$$
so $f\in \bF_{g,\mathbf{\infty}} (R,Q)$ simply means that the common density of observations $\mathfrak{p}$ is uniformly bounded by $Q$.
It is easy to see that if $\alpha=1$ and $\|g\|_\infty<\infty$ then $\bF_{g,\mathbf{\infty}} (R,Q)=\bF_g(R)$ for any $Q\geq R\|g\|_\infty$.

Let us briefly discuss another example. Let $r>1$ and $L<\infty$ be arbitrary but a priory chosen numbers. Assume that the considered collection of anisotropic Nikol'skii classes obeys the following restrictions: $\vec{r}\in [r,\infty]^d$ and $\vec{L}\in (0,L]^d$. Suppose also that $\|g\|_s<\infty$, where $1/s=1-1/r$. Then, there exists $Q_0$ completely determined by $r,L$ and $R$ such that $\bN_{\vec{r},d}\big(\vec{\beta},\vec{L}\big)\cap\bF_{g,\mathbf{\infty}} (R,Q)=\bN_{\vec{r},d}\big(\vec{\beta},\vec{L}\big)\cap\bF_{g} (R)$ for any $Q>Q_0\|g\|_s$.

\vskip0.1cm

Additionally, we will study the adaptive estimation over the collection
 $$
 \bF_\vartheta=\bN_{\vec{r},d}\big(\vec{\beta},\vec{L}\big)\cap\bF_g(R)\cap\bB_{\infty,d}(Q),\; \vartheta=\big(\vec{\beta},\vec{r},\vec{L},R,Q\big).
 $$
We will show that the boundedness
of the underlying function allows to improve considerably
 the accuracy of estimation.

\subsection{Historical notes}

The minimax adaptive estimation is a very active area of mathematical statistics and  the interested reader can find a very detailed overview as well as several open problems in adaptive estimation in the recent paper,
\cite{lepski15}. Below we will discuss only the articles whose results are relevant to our consideration, i.e. the density setting under $\bL_p$-loss, from a minimax perspective.

\smallskip

Let us  start with the following remark. If one assumes additionally that $f,g\in\mP\big(\bR^d\big)$ the convolution structure density model can be interpreted as follows. The observations $Z_i\in\bR^d, i=1,\ldots, n,$ can be written as a sum of two independent random vectors, that is,
\begin{equation}
\label{eq:observation-scheme}
Z_i=X_i+\epsilon_iY_i,\quad i=1,\ldots,n,
\end{equation}
where $X_i, i=1,\ldots,n,$ are {\it i.i.d.} $d$-dimensional  random vectors with common density $f$ to be estimated.
The noise variables  $Y_i, i=1,\ldots,n,$
are {\it i.i.d.} $d$-dimensional random vectors with known common density $g$. At last $\e_i\in\{0,1\}, i=1,\ldots,n,$
are {\it i.i.d.} Bernoulli random variables with $\bP(\e_1=1)=\alpha$, where $\alpha\in [0,1]$ is supposed to be known.
The sequences  $\{X_i, i=1,\ldots,n\}$, $\{Y_i, i=1,\ldots,n\}$ and $\{\epsilon_i, i=1,\ldots,n\}$ are supposed to be mutually independent.

The  observation scheme (\ref{eq:observation-scheme}) can be viewed as the generalization of two classical statistical models. Indeed, the case $\alpha=1$ corresponds to the standard deconvolution model $Z_i=X_i+Y_i,\; i=1,\ldots,n$. Another "extreme" case $\alpha=0$ correspond to the direct observation scheme $Z_i=X_i,\; i=1,\ldots,n$.
The "intermediate" case $\alpha\in (0,1)$, considered for the first time in \cite{hesse}, is understood as partially contaminated observations.



\paragraph{\textsf{Direct case,} $\alpha=0$}

There is a vast literature dealing with minimax and minimax adaptive density estimation, see for example, \cite{Efr}, \cite{Has-Ibr}, \cite{golubev92}, \cite{Donoho},  \cite{dev-lug97}, \cite{rigollet}, \cite{rigollet-tsybakov}, \cite{samarov}, \cite{birge}, \cite{GN-1}, \cite{akakpo}, \cite{Gach}, \cite{lepski13a},
among many others. Special attention was paid to the estimation of densities with unbounded support, see \cite{Juditsky}, \cite{patricia}.
The most developed results  can be found in
\cite{GL11}, \cite{GL14} and  in Section \ref{sec:adaptive-results-deconv} we will compare in detail our results with those obtained in these papers.

\paragraph{\textsf{Intermediate case,} $\alpha\in (0,1)$}

To the best of our knowledge, adaptive estimation in the case of partially contaminated observations has not been studied yet. We were able to find only two papers dealing with minimax estimation. The first one is \cite{hesse}
(where the discussed model was introduced in dimension $1$) in which the author evaluated the $\bL_\infty$-risk of the proposed estimator over a functional class formally corresponding to the Nikol'skii class $\bN_{\infty,1}(2,1)$. In \cite{Yuana-Chenb} the latter result was developed to the multidimensional setting, i.e. to the minimax estimation on  $\bN_{\infty,d}\big(\vec{2},1\big)$.  The most intriguing fact is
that the accuracy of estimation in partially contaminated noise is the same as in the direct observation scheme.
However none of these articles  studied the optimality of the proposed estimators.
 We will come back to the aforementioned papers in Section \ref{sec:subsec-assumptions on the distribution of the noise-from-LW} in order to compare the assumptions imposed on the noise density $g$.

\paragraph{\textsf{Deconvolution case,} $\alpha=1$}

First let us remark that the behavior of the Fourier transform of the  function $g$ plays an important role in all the works dealing with deconvolution. Indeed ill-posed problems correspond to Fourier transforms decaying towards zero. Our results will be established for "moderately" ill posed problems, so we detail only results in papers studying that type of operators. This assumption means that there exist $\vec{\mu}=(\mu_1,\ldots,\mu_d)\in (0,\infty)^{d}$ and $\Upsilon_1>0, \Upsilon_2>0$ such that the Fourier transform $\check{g}$ of $g$ satisfies:
\begin{eqnarray}
\label{eq:illposed-moderate}
&&\Upsilon_{1} \prod_{j=1}^d(1+\ttt^2_j)^{-\frac{\mu_j}{2}} \leq \big|\check{g}(\ttt)\big|\leq \Upsilon_{2} \prod_{j=1}^d(1+\ttt^2_j)^{-\frac{\mu_j}{2}},\quad\forall \ttt=(\ttt_1,\ldots,\ttt_d)\in\bR^d.
\end{eqnarray}

Some minimax and minimax adaptive results in dimension 1 over different classes of smooth functions can be found in particular in  \cite{sc90}, \cite{f91},
\cite{f93},
 \cite{pv99}, \cite{fk02},  \cite{c06},  \cite{hm07}, \cite{m09},
\cite{Lounici-Nickl}, \cite{kerk11}.

There are very few results in the multidimensional setting. It seems that \cite{Masry93} was the first paper where the deconvolution problem was studied for multivariate densities. It is worth noting that \cite{Masry93} considered more general weakly dependent observations and this paper formally does not deal with the minimax setting. However the results obtained in this paper could be formally compared with the estimation under $\bL_\infty$-loss over the \textsf{isotropic} H\"older class of regularity $2$, i.e. $\bN_{\infty,d}\big(\vec{2},1\big)$ which is exactly the same setting as in \cite{Yuana-Chenb} in the case of partially contaminated observations.  Let us also remark that there is no lower bound result in \cite{Masry93}. The most developed results in the deconvolution model were obtained in \cite{comte} and \cite{rebelles16}
 and  in Section \ref{sec:adaptive-results-deconv} we will compare in detail our results with those obtained in these papers.

\subsection{Lower bound for the minimax $\bL_p$-risk}
\label{sec:subsec-Lower bound for minimax-risk-deconv}
We have seen that the problem of optimal adaptation over the collection $\big\{\bF_\vartheta,\vartheta\in\Theta\big\}$ is formulated as the "attainability" of the family of minimax risks $\big\{\phi_n(\bF_\vartheta), \vartheta\in\Theta\big\}$ by a single estimator.
Although it is not necessary, the following "two-stage" approach is used for the majority of problems related to the minimax adaptive  estimation.
The first step consists in finding a lower bound for $\phi_n(\bF_\vartheta)$ for any $\vartheta\in\Theta$
while the second one consists in constructing an estimator "attaining", at least asymptotically, this bound.
We adopt this strategy in our investigations and below we present several lower bound results recently obtained in \cite{LW15}.

\subsubsection{\textsf{Assumptions on the function $g$ imposed in \cite{LW15}}}
\label{sec:subsec-assumptions on the distribution of the noise-from-LW}
Let $\mJ^*$ denote the set of all subsets of  $\{1,\ldots,d\}$. Set $\mJ=\mJ^*\cup\emptyset$ and for any $J\in\mJ$ let $|J|$ denote the cardinality of $J$
while $\{j_1<\cdots< j_{|J|}\}$ denotes its elements.

 For any $J\in\mJ^*$ define the operator
$
\mD^{J}=\frac{\partial^{|J|}}{\partial\ttt_{j_1}\cdots\partial\ttt_{j_{|J|}}}
$
and let  $\mD^{\emptyset}$ denote the identity operator. For any $I,J\in\mJ$ define $\mD^{I,J}=\mD^{I}\big(\mD^J\big)$ and note that obviously $\mD^{I,J}=\mD^{J,I}$.

\begin{assumption}[$\alpha\neq 1$]

\label{ass2:ass-on-noise-lower-bound}

$\mD^{J}\check{g}$ exists for any $J\in\mJ^*$ and
$
\;\sup_{J\in\mJ^*}\big\|\mD^{J}\check{g}\big\|_\infty<\infty;
$

\end{assumption}

\begin{assumption}[$\alpha=1$]
\label{ass1:ass-on-noise-lower-bound}
 $\mD^{J}\check{g}$ exists for any $J\in\mJ^*$ and
$
\sup_{J\in\mJ^*}\big\|\check{g}^{-1}\mD^{J}\check{g}\big\|_\infty<\infty.
$
Moreover,  there exists $\vec{\mu}=(\mu_1,\ldots,\mu_d)\in (0,\infty)^{d}$ and $\Upsilon>0$  such that

\vskip0.2cm

\centerline{$
|\check{g}(\ttt)|\leq \Upsilon\prod_{j=1}^d(1+\ttt^2_j)^{-\frac{\mu_j}{2}},\quad\forall \ttt=(\ttt_1,\ldots,\ttt_d)\in\bR^d.
$}
\end{assumption}

\begin{assumption}[$\alpha=1$]
\label{ass4:ass-on-noise-lower-bound}
$g$ is a bounded function.
\end{assumption}

\begin{assumption}[$\alpha=1$]
\label{ass3:ass-on-noise-lower-bound}

$\mD^{I,J}\check{g}$ exists for any $I,J\in\mJ$ and $\;\sup_{I,J\in\mJ}\big\|\mD^{I,J}\big(\check{g}\big)\big\|_1<\infty$. Moreover
\begin{eqnarray*}
\sup_{J\in\mJ^*}\int_{\bR^d}g(z)\Big(\prod_{j\in J}z^{2}_j\Big)\rd z<\infty.
\end{eqnarray*}

\end{assumption}

It is worth noting that all the bounds in \cite{LW15} are obtained under Assumptions \ref{ass2:ass-on-noise-lower-bound} and \ref{ass1:ass-on-noise-lower-bound}. Assumption \ref{ass4:ass-on-noise-lower-bound} is used when the estimation of unbounded functions is considered; we come back to this assumption in Section  \ref{sec:subsubsec:unbounded case}.


As to  Assumption  \ref{ass3:ass-on-noise-lower-bound}, it seems purely technical and does not appear in upper bound results. We also recall that the lower bounds in \cite{LW15} are proved under the condition: $g\in\mP\big(\bR^d\big)$.

\subsubsection{\textsf{Some lower bounds from \cite{LW15}}}

Set $\vec{\blg}(\alpha)=\vec{\mu}$, $\alpha=1$, $\vec{\blg}(\alpha)=(0,\ldots,0)$, $\alpha\in [0,1)$, and introduce for any
$\vec{\beta}\in (0,\infty)^d$, $\vec{r}\in [1,\infty]^d$ and $\vec{L}\in (0,\infty)^d$
 the following quantities.
\begin{eqnarray}
 \label{eq:beta&omega}
\frac{1}{\beta(\alpha)}=\sum_{j=1}^d\frac{2\blg_j(\alpha)+1}{\beta_j},\quad\frac{1}{\omega(\alpha)}=\sum_{j=1}^d\frac{2\blg_j(\alpha)+1}{\beta_j r_j},
\quad L(\alpha)=\prod_{j=1}^dL_j^{\frac{2\blg_j(\alpha)+1}{\beta_j}}.
\end{eqnarray}
Define for any $1\leq s\leq\infty$ and $\alpha\in[0,1]$
\begin{eqnarray}
 \label{eq:kappa&tau}
 \kappa_\alpha(s)=\omega(\alpha)(2+1/\beta(\alpha))-s\qquad  \tau(s)=1-1/\omega(0)+1/(s\beta(0)).
\end{eqnarray}

\noindent{\textsf{General case.}
Remind that  $
z(\alpha)=\omega(\alpha)
(2+1/\beta(\alpha))\beta(0)\tau(\infty)+1
$,   $p^*=\big[\max_{l=1,\ldots,d}r_l\big]\vee p$. Set
\begin{eqnarray}
\label{eq:varrho(alpha)}
 \varrho(\alpha) &=&\left\{\begin{array}{cclcc}
    \frac{1-1/p}{1-1/\omega(\alpha)+1/\beta(\alpha)},\; & \kappa_\alpha(p)> p\omega(\alpha) ;
    \\*[2mm]
\frac{\beta(\alpha)}{2\beta(\alpha)+1},\; &   0<\kappa_\alpha(p)\leq p\omega(\alpha);
\\*[2mm]
\frac{\tau(p)\omega(\alpha)\beta(0)}{z(\alpha)},\,  &\kappa_\alpha(p)\leq 0,\;\; \tau(p^*)> 0;
\\*[2mm]
\frac{\omega(\alpha)(1-p^*/p)}{\kappa_\alpha(p^*)},\; & \kappa_\alpha(p)\leq 0,\;\tau(p^*)\leq 0.
            \end{array}
\right.
\end{eqnarray}
Here and later we assume $0/0=0$, which implies in particular that $\frac{\omega(\alpha)(1-p^*/p)}{\kappa_\alpha(p^*)}=0$ if $p^*=p$ and $\kappa_\alpha(p)=0$. Recall also that $\kappa_\alpha(p^*)/p^*=-1$ if $p^*=\infty$.
Put at last
\begin{eqnarray*}
\bld_n &=& \left\{
\begin{array}{ccl}
L(\alpha)n^{-1},\quad & \kappa_\alpha(p)> 0 ;
\\*[2mm]
L(\alpha)n^{-1}\ln (n),\quad &\kappa_\alpha(p)\leq 0,\;\;\tau(p^*)\leq 0;
\\*[2mm]
[L(0)]^{-\frac{\kappa_\alpha(p)}{\omega(\alpha)p\tau(p)}}L(\alpha)n^{-1}\ln (n),\quad &\kappa_\alpha(p)\leq 0,\;\;\tau(p^*)> 0.
\end{array}
\right.
\end{eqnarray*}

\begin{theorem}[\cite{LW15}]
\label{th1:lower-bound-deconvolution}
Let   $L_0>0$ and $1\leq p<\infty$ be fixed.

\vskip0.2cm

Then for any  $\vec{\beta}\in (0,\infty)^d,\; \vec{r}\in [1,\infty]^d$,  $\vec{L}\in [L_0,\infty)^d$, $\vec{\mu}\in (0,\infty)^d$, $R>1$ and
$g\in\mP\big(\bR^d\big)$, satisfying Assumptions \ref{ass2:ass-on-noise-lower-bound}--\ref{ass3:ass-on-noise-lower-bound},
there exists $c>0$ independent of $\vec{L}$ such that

$$
\liminf_{n\to \infty}\;\inf_{\tilde{f}_n}\sup_{f\in\bN_{\vec{r},d}\big(\vec{\beta},\vec{L}\big)\cap\bF_g(R)}\bld_n^{-\varrho(\alpha)}\cR^{(n)}_p\big[\tilde{f}_n; f\big]\geq c,
$$
where the infimum is taken over all possible estimators.

\end{theorem}

Following the terminology used in \cite{LW15}, we will call the set of parameters satisfying $\kappa_\alpha(p)> p\omega(\alpha)$ the  \textsf{tail zone},  satisfying $0<\kappa_\alpha(p)\leq p\omega(\alpha)$ the  \textsf{dense zone} and satisfying $\kappa_\alpha(p)\leq 0$ the \textsf{sparse zone}. In its turn, the latter zone is divided into two sub-domains: the \textsf{sparse zone 1} corresponding to
$\tau(p^*)> 0$ and the \textsf{sparse zone 2} corresponding to
$\tau(p^*)\leq 0$.

\smallskip

\noindent\textsf{Bounded case.}
Introduce
\begin{eqnarray}
\label{eq:rho(alpha)}
 \rho(\alpha) &=&\left\{\begin{array}{clc}
    \frac{1-1/p}{1-1/\omega(\alpha)+1/\beta(\alpha)},\; & \kappa_\alpha(p)> p\omega(\alpha) ;
    \\*[2mm]
\frac{\beta(\alpha)}{2\beta(\alpha)+1},\; &   0<\kappa_\alpha(p)\leq p\omega(\alpha);
\\*[2mm]
\frac{\tau(p)\omega(\alpha)\beta(0)}{z(\alpha)},\;  &\kappa_\alpha(p)\leq 0,\;\; \tau(\infty)> 0;
\\*[2mm]
\frac{\omega(\alpha)}{p},\; & \kappa_\alpha(p)\leq 0,\;\tau(\infty)\leq 0.
            \end{array}
\right.
\end{eqnarray}

\begin{theorem}[\cite{LW15}]
\label{th2:lower-bound-deconvolution}
Let   $L_0>0$ and $1\leq p<\infty$ be fixed.

\vskip0.1cm

Then for any  $\vec{\beta}\in (0,\infty)^d,\; \vec{r}\in [1,\infty]^d$,  $\vec{L}\in [L_0,\infty)^d$, $Q>0$, $\vec{\mu}\in (0,\infty)^d$ and
$g\in\mP\big(\bR^d\big)$, satisfying Assumptions \ref{ass2:ass-on-noise-lower-bound} and \ref{ass1:ass-on-noise-lower-bound}
there exists $c>0$ independent of $\vec{L}$ such that
$$
\liminf_{n\to \infty}\;\inf_{\tilde{f}_n}\sup_{f\in\bN_{\vec{r},d}\big(\vec{\beta},\vec{L}\big)\cap\mP\big(\bR^d\big)\cap\bB_{\infty,d}(Q)}
\bld_n^{-\rho(\alpha)}\cR^{(n)}_p\big[\tilde{f}_n; f\big]\geq c,
$$
where the infimum is taken over all possible estimators.
\end{theorem}

\subsection{Assumptions on the function $g$}

The selection rule from the family of linear estimators, the $\bL_p$-norm oracle inequalities obtained in Part I and all the adaptive results presented in the paper are established under the following condition imposed on the function $g$.

\begin{assumption}
\label{ass1:ass-on-noise-upper-bound}

\smallskip

(1) if $\alpha\neq 1$ then there exists  $\e>0$ such that

\vskip0.2cm

\centerline{$
\big|1 - \alpha +\alpha\check{g}(t)\big|\geq \e, \quad \forall t\in\bR^d;
$}

\vskip0.2cm

(2) if $\alpha=1$ then there exists  $\vec{\mu}=(\mu_1,\ldots,\mu_d)\in (0,\infty)^{d}$ and $\Upsilon_0>0$  such that
\vskip0.2cm

\centerline{$
|\check{g}(t)|\geq \Upsilon_0\prod_{j=1}^d(1+\ttt^2_j)^{-\frac{\mu_j}{2}},\quad\forall \ttt=(\ttt_1,\ldots,\ttt_d)\in\bR^d.
$}

\vskip0.2cm

\end{assumption}

Comparing this condition with Assumption \ref{ass1:ass-on-noise-lower-bound} from Section \ref{sec:subsec-assumptions on the distribution of the noise-from-LW}, we can assert that both
are coherent  if $\alpha=1$.
Indeed, in this case, we come the following assumption, which is well-known in the literature:
\vskip0.2cm

\centerline{$
 \Upsilon_0\prod_{j=1}^d(1+t^2_j)^{-\frac{\mu_j}{2}}\leq |\check{g}(t)|\leq \Upsilon\prod_{j=1}^d(1+\ttt^2_j)^{-\frac{\mu_j}{2}},\quad\forall t\in\bR^d.
$}

\vskip0.2cm

\noindent referred to as a {\it moderately ill-posed} statistical problem, cf. (\ref{eq:illposed-moderate}). In particular, the assumption is checked for the centered multivariate Laplace law.

Note first that Assumption \ref{ass1:ass-on-noise-upper-bound} is in some sense weaker than Assumption \ref{ass2:ass-on-noise-lower-bound} when $\alpha\in(0,1)$, since it does not require regularity properties of the function $g$. Moreover both assumptions are not too restrictive. They are verified for many distributions, including centered  multivariate Laplace and Gaussian ones.
Note also that Assumption \ref{ass1:ass-on-noise-upper-bound} always holds with $\e=1-2\alpha$ if $\alpha<1/2$.
Additionally, it holds with $\e=1-\alpha$ if $\check{g}$ is a real positive function. The latter is true, in particular, for any probability law obtained by an even number
of convolutions of a symmetric distribution with itself.

Next,  our Assumption \ref{ass1:ass-on-noise-upper-bound} is weaker than the conditions imposed in \cite{hesse} and  \cite{Yuana-Chenb}. In these papers $\check{g}\in\bC^{(2)}\big(\bR^d\big)$,
$\check{g}(t)\neq 0$ for any $t\in\bR^d$ and
$$
\big|1 - \alpha +\alpha\check{g}(t)\big|\geq 1-\alpha, \quad \forall t\in\bR^d.
$$

\section{Adaptive estimation over the scale of anisotropic Nikol'skii classes}
\label{sec:adaptive-results-deconv}

We start this section by recalling the definition of the pointwise selection rule proposed in Part I.

\subsection{Pointwise selection rule}

Let $K:\bR^d\to\bR$ be a continuous  function belonging to $\bL_1\big(\bR^d\big)$ such that $\int_{\bR}K=1$.
Set  $\cH=\big\{e^{k},\; k\in\bZ\big\}$ and let
$
\cH^d=\big\{\vec{h}=(h_1,\ldots,h_d):\; h_j\in\cH, j=1,\ldots,d\big\}.
$
Recall that $
\cH^d_{\text{isotr}}=\big\{\vec{h}\in\cH^d:\; \vec{h}=(h,\ldots,h), \; h\in\cH\big\}.
$
Set $ V_{\vec{h}}=\prod_{j=1}^dh_j$ and let for any $\vec{h}\in\cH^d$
$$
K_{\vec{h}}(t)=V^{-1}_{\vec{h}}K\big(t_1/h_1,\ldots,t_d/h_d\big),\; t\in\bR^d.
$$
Later on for any $u,v\in\bR^d$ the operations and relations $u/v$, $uv$, $u\vee v$,$u\wedge v$,
$u\geq v$, $au, a\in\bR,$ are understood in coordinate-wise sense. In particular $u\geq v$ means that $u_j\geq v_j$ for any $j=1,\ldots, d$.

For any $\vec{h}\in (0,\infty)^d$ let  $M\big(\cdot,\vec{h}\big)$
satisfy the operator equation
\begin{eqnarray}
\label{eq:def-kernel-M}
K_{\vec{h}}(y)=(1-\alpha)M\big(y,\vec{h}\big)+\alpha\int_{\bR^d}g(t-y)M\big(t,\vec{h}\big)\rd t,\quad  y\in\bR^d.&&
\end{eqnarray}
Introduce for any $\vec{\mathrm{h}}\in\cH^d$ and $x\in\bR^d$
\begin{gather*}
\widehat{f}_{\vec{\mathrm{h}}}(x)=n^{-1}\sum_{i=1}^n M\big(Z_i-x,\vec{\mathrm{h}}\big), \qquad \widehat{\sigma}^2\big(x,\vec{\mathrm{h}}\big)=\frac{1}{n}\sum_{i=1}^nM^2\big(Z_i-x,\vec{\mathrm{h}}\big);
\\
\widehat{U}_n\big(x,\vec{\mathrm{h}}\big)=\sqrt{\frac{2\lambda_n\big(\vec{\mathrm{h}}\big)\widehat{\sigma}^2\big(x,\vec{\mathrm{h}}\big)}{n}}+\frac{4 M_\infty\lambda_n\big(\vec{\mathrm{h}}\big)}{3n\prod_{j=1}^d\mathrm{h}_j(\mathrm{h}_j\wedge 1)^{\blg_j(\alpha)}},
\end{gather*}
where $M_\infty=\big[(2\pi)^{-d}\big\{\e^{-1}\big\|\check{K}\big\|_1\mathrm{1}_{\alpha\neq 1}+\Upsilon_0^{-1}\mathbf{k}_1\mathrm{1}_{\alpha=1}\big\}\big]\vee 1$ and

\vskip-0.5cm

\begin{gather*}
\lambda_n\big(\vec{\mathrm{h}}\big)=4\ln(M_\infty)+6\ln{(n)}+(8p+26)\sum_{j=1}^d\big[1+\blg_j(\alpha)\big]\big|\ln(\mathrm{h}_j)\big|.
\end{gather*}

\vskip-0.2cm

 Let $\bH$ be an arbitrary subset of $\cH^d$. For any $\vec{h}\in\bH$ and $x\in\bR^d$ introduce
\begin{eqnarray}
\label{eq:pointwise-rule-deconvolution}
&&\widehat{\cR}_{\vec{h}}(x)=\sup_{\vec{\eta}\in\bH}\Big[\big|\widehat{f}_{\vec{h}\vee\vec{\eta}}(x)-\widehat{f}_{\vec{\eta}}(x)\big|
-4\widehat{U}_n\big(x,\vec{h}\vee\vec{\eta}\big)-4\widehat{U}_n\big(x,\vec{\eta}\big)\Big]_+;
\\
&&\widehat{U}^*_n\big(x,\vec{h}\big)=\sup_{\vec{\eta}\in\bH:\; \vec{\eta}\geq\vec{h}}\widehat{U}_n\big(x,\vec{\eta}\big),
\end{eqnarray}
and define
$\vec{\mathbf{h}}(x)=\arg\inf_{\vec{h}\in\bH}\Big[\widehat{\cR}_{\vec{h}}(x)+8\widehat{U}^*_n\big(x,\vec{h}\big)\Big].$

Our final estimator is $\widehat{f}_{\vec{\mathbf{h}}(x)}(x),\;x\in\bR^d$ and we will call (\ref{eq:pointwise-rule-deconvolution}) the {\it pointwise selection rule}.

\begin{remark}
Note that the estimator $\widehat{f}_{\vec{\mathbf{h}}}$ depends on $\bH$ and
later on we will consider two choices of the parameter set $\bH$, namely $\bH=\cH^d$ and $\bH=\cH^d_{\text{isotr}}$. So, to present our results we will write $\widehat{f}_{\vec{\mathbf{h}},\bH}$ in order to underline the aforementioned dependence.
The choice $\bH=\cH^d$  will be used when the adaptation is studied over anisotropic Nikol'skii classes while $\bH=\cH^d_{\text{isotr}}$ will be used when the considered scale consists of isotropic classes.

\end{remark}

\subsection{Anisotropic Nikol'skii classes}
\label{sec:nikolski}
Let $(\mathbf{e}_1,\ldots,\mathbf{e}_d)$ denote the canonical basis of $\bR^d$.
 For some function $G:\bR^d\to \bR^1$ and
real number $u\in \bR$ define
{\em the first order difference operator with step size $u$ in direction of the variable
$x_j$} by
\[
 \Delta_{u,j}G (x)=G(x+u\mathbf{e}_j)-G(x),\;\;\;j=1,\ldots,d.
\]
By induction,
the $k$-th order difference operator with step size $u$ in direction of the variable $x_j$ is
defined as
\begin{equation}
\label{eq:Delta}
 \Delta_{u,j}^kG(x)= \Delta_{u,j} \Delta_{u,j}^{k-1} G(x) = \sum_{l=1}^k (-1)^{l+k}\binom{k}{l}\Delta_{ul,j}G(x).
\end{equation}
\begin{definition}
\label{def:nikolskii}
For given  vectors $\vec{r}=(r_1,\ldots,r_d)\in [1,\infty]^d$ $\vec{\beta}=(\beta_1,\ldots,\beta_d)\in (0,\infty)^d$
 and $\vec{L}=(L_1,\ldots, L_d)\in (0,\infty)^d$ we
say that a function $G:\bR^d\to \bR^1$ belongs to the anisotropic
Nikolskii class $\bN_{\vec{r},d}\big(\vec{\beta},\vec{L}\big)$ if
\vskip0.1cm

{\rm (i)}\;\; $\|G\|_{r_j}\leq L_{j}$ for all $j=1,\ldots,d$;

\vskip0.1cm

{\rm (ii)}\;\;
for every $j=1,\ldots,d$ there exists natural number  $k_j>\beta_j$ such that
\begin{equation*}
\label{eq:Nikolski}
 \Big\|\Delta_{u,j}^{k_j} G\Big\|_{r_j} \leq L_j |u|^{\beta_j},\;\;\;\;
\forall u\in \bR,\;\;\;\forall j=1,\ldots, d.
\end{equation*}
\end{definition}

If $\beta_j=\blb\in (0,\infty), r_j=\mathbf{r}\in [1,\infty]$ and $L_j=\mathbf{L}\in (0,\infty)$ for any $j=1,\ldots, d$
the corresponding Nikolskii class, denoted furthermore $\bN_{\mathbf{r},d}(\blb,\mathbf{L})$, is called isotropic.

\subsection{Construction of kernel $K$} First, we recall that all results concerning the $\bL_p$ risk of the pointwise selection rule, established in Part I,
are proved under the following assumption imposed on the kernel $K$.

\begin{assumption}
\label{ass:on-kernel-deconvolution}
There exist $\mathbf{k}_1>0$ and $\mathbf{k}_2>0$ such that
\begin{eqnarray}
\int_{\bR^d}\big|\check{K}(t)\big|\prod_{j=1}^d(1+t^2_j)^{\frac{\blg_j(\alpha)}{2}}\rd t\leq \mathbf{k}_1,\quad
\int_{\bR^d}\big|\check{K}(t)\big|^2\prod_{j=1}^d(1+t^2_j)^{\blg_j(\alpha)}\rd t\leq \mathbf{k}^2_2.&&
\end{eqnarray}

\end{assumption}

Next, we will use the following specific kernel $K$ in the definition of the estimator's family $\big\{\widehat{f}_{\vec{\mathrm{h}}}(\cdot),\;\vec{\mathrm{h}}\in\cH^d \big\}$
[see, e.g., \cite{lepski-kerk} or \cite{GL14}].

\vskip0.1cm

 Let  $\ell$ be an integer number,
and let $\cK:\bR^1\to \bR^1$ be a \textsf{compactly supported} continuous function satisfying  $\int_{\bR^1} \cK(y)\rd y=1$,
and $\cK\in\bC(\bR^1)$. Put
\begin{equation}
\label{eq:w-function}
 \cK_\ell(y)=\sum_{i=1}^\ell \binom{\ell}{i} (-1)^{i+1}\frac{1}{i}\cK\Big(\frac{y}{i}\Big),
\end{equation}
and add the following structural condition to Assumption \ref{ass:on-kernel-deconvolution}.
\begin{assumption}
\label{ass2:on-kernel-deconvolution}
$
K(x)=\prod_{j=1}^d \cK_\ell(x_j),\; \forall x\in\bR^d.
$
\end{assumption}

The kernel $K$ constructed in this way is bounded, compactly  supported, belongs to $\bC(\bR^d)\cap\bL_1(\bR^d)$
and satisfies $\int_{\bR^d} K=1$. Some examples of  kernels satisfying simultaneously Assumptions \ref{ass:on-kernel-deconvolution} and \ref{ass2:on-kernel-deconvolution} can be found for instance  in \cite{comte}.

\subsection{Main results}
\label{sec:adap-nikol-deconv}


Introduce the following notations: $\delta_n=L(\alpha)n^{-1}\ln (n)$ and
\begin{eqnarray*}
t(\bH)=\left\{
\begin{array}{cc}
d-1,\;& \bH=\cH^d;
\\*[2mm]
0,\; &\bH=\cH^d_{\text{isotr}},
\end{array}
\right.
\quad
\mb_n(\bH) &=& \left\{
\begin{array}{ccl}
[\ln(n)]^{t(\bH)},\quad &\qquad\kappa_\alpha(p)>p\omega(\alpha);
\\*[2mm]
\ln^{\frac{1}{p}} (n)\vee[\ln(n)]^{t(\bH)},\quad &\qquad\kappa_\alpha(p)=p\omega(\alpha);
\\*[2mm]
\ln^{\frac{1}{p}} (n),& \kappa_\alpha(p)=0;
\\*[2mm]
1,& \text{otherwise},
\end{array}
\right.
\end{eqnarray*}

\subsubsection{\textsf{Bounded case}} The first  problem we address  is the adaptive estimation over the collection of the functional classes
$\big\{\bN_{\vec{r},d}\big(\vec{\beta},\vec{L}\big)\cap\bF_g(R)\cap\bB_{\infty,d}(Q)\big\}_{\vec{\beta},\vec{r},\vec{L},R,Q}.\;$

\vskip0.1cm

As it was conjectured  in \cite{LW15},
the boundedness of the function belonging to $\bN_{\vec{r},d}\big(\vec{\beta},\vec{L}\big)\cap\bF_g(R)$ is a minimal condition allowing to eliminate the inconsistency zone. The results obtained in Theorem \ref{th:adaptive-L_p-bounded-deconv} below together with those from Theorem
\ref{th2:lower-bound-deconvolution} confirm this conjecture.

\begin{theorem}
\label{th:adaptive-L_p-bounded-deconv}
Let  $\alpha\in [0,1]$,   $\ell\in\bN^*$ and $g\in\bL_1\big(\bR^d\big)$, satisfying Assumption \ref{ass1:ass-on-noise-upper-bound}, be fixed. Let $K$ satisfy Assumptions \ref{ass:on-kernel-deconvolution} and \ref{ass2:on-kernel-deconvolution}.

1) Then for any $p\in(1,\infty)$, $Q>0$, $R>0$, $ L_0>0$, $\vec{\beta}\in (0,\ell]^d$, $\vec{r}\in (1,\infty]^d$ and $\vec{L}\in [L_0,\infty)^d$
there exists $C<\infty$, independent of $\vec{L}$,  such that:
$$
\limsup_{n\to \infty}\sup_{f\in\bN_{\vec{r},d}\big(\vec{\beta},\vec{L}\big)\cap\bF_g(R)\cap\bB_{\infty,d}(Q)}\mb_n\big(\cH^d\big)^{-1}\delta_n^{-\rho(\alpha)}
\cR^{(n)}_p\big[\widehat{f}_{\vec{\mathbf{h}},\bH}; f\big]\leq C,
$$
where  $\rho(\alpha)$ is defined in (\ref{eq:rho(alpha)}).

2) For any $p\in(1,\infty)$, $Q>0$, $R>0, L_0>0$, $\blb\in (0,\ell]$, $\mathbf{r}\in [1,\infty]$ and $\mathbf{L}\in [L_0,\infty)$
there exists $C<\infty$, independent of $\mathbf{L}$,  such that:
$$
\limsup_{n\to \infty}\sup_{f\in\bN_{\mathbf{r},d}\big(\blb,\mathbf{L}\big)\cap\bF_g(R)\cap\bB_{\infty,d}(Q)}\mb_n\big(\cH^d_{\text{isotr}}\big)^{-1}\delta_n^{-\rho(\alpha)}
\cR^{(n)}_p\big[\widehat{f}_{\vec{\mathbf{h}},\cH^d_{\text{isotr}}}; f\big]\leq C,
$$
\end{theorem}

Some remarks are in order.
$\mathbf{1^0}.\;$ Our estimation procedure is completely data-driven, i.e. independent of $\vec{\beta},\vec{r}, \vec{L}, R$, $Q$,
and the assertions of Theorem \ref{th:adaptive-L_p-bounded-deconv} are completely new if $\alpha\neq 0$.
Comparing the results obtained in Theorems \ref{th2:lower-bound-deconvolution} and \ref{th:adaptive-L_p-bounded-deconv} we can assert that our estimator is optimally-adaptive if $\kappa_\alpha(p)<0$ and nearly optimally adaptive if $0<\kappa_\alpha(p)<p\omega(\alpha)$. The construction of an estimation procedure which would be optimally-adaptive when $\kappa_\alpha(p)\geq 0$ is an open problem, and we conjecture that the lower bounds for the asymptotics of the minimax risk found in Theorem \ref{th2:lower-bound-deconvolution} are sharp in order. This conjecture in the case $\alpha=1$ is partially confirmed by the results obtained in \cite{comte} and \cite{rebelles16}. Since both articles deal with the estimation of unbounded functions we will discuss them in the next section.

It is worth noting that all the previous statements are true not only for the convolution structure density model but also, in view of Theorem \ref{th2:lower-bound-deconvolution}, for the observation scheme (\ref{eq:observation-scheme}) as well.

\smallskip

$\mathbf{2^0}.\;$ We note that the asymptotic of the minimax risk under  partially contaminated observations, $\alpha\in (0,1)$, is independent of $\alpha$ and coincides with the asymptotic of the risk in the direct observation model, $\alpha=0$. For the first time this phenomenon was discovered in \cite{hesse} and \cite{Yuana-Chenb}.  In the very recent paper \cite{lepski16}, in the particular case  $\vec{r}=(p,\ldots,p)$, $p\in (1,\infty)$ the optimally adaptive estimator
 was built.  It is easy to check that independently of the value of $\vec{\beta}$ and $\vec{\mu}$, the corresponding set of parameters belongs to the \textsf{dense} zone.  Note however that our estimator is only optimally-adaptive in this zone, but it is applied to a much more general collection of functional classes. It is worth noting that the estimator procedure, used  in \cite{lepski16}, has nothing in common with our pointwise selection rule.

 \smallskip

$\mathbf{3^0}.\;$
 As to the direct observation scheme, $\alpha=0$,  our results coincide with those obtained recently in \cite{GL14}, when
 $p\omega(0)>\kappa_0(p)$.
  However,  for the tail zone $p\omega(0)\leq\kappa_0(p)$, our bound is slightly better since the bound obtained in the latter paper contains an additional factor $\ln^{\frac{d}{p}}(n)$. It is interesting to note that although both estimator constructions are based upon local selections from the family of kernel estimators, the selection rules are different.

\smallskip

$\mathbf{4^0}.\;$ Let us finally discuss the results corresponding
to the \textsf{tail zone},
$\kappa_\alpha(p)>p\omega(\alpha)$.
First, the lower bound for the minimax risk is given by $[L(\alpha)n^{-1}]^{\rho(\alpha)}$
 while the accuracy provided by our estimator is
$$
\ln^{\frac{d-1}{p}}(n)[L(\alpha)n^{-1}\ln (n)]^{\rho(\alpha)}.
$$

As we mentioned above, the passage from $[L(\alpha)n^{-1}]^{\rho(\alpha)}$ to
$[L(\alpha)n^{-1}\ln (n)]^{\rho(\alpha)}$ seems to be  an unavoidable payment  for the application of a local selection scheme. It is interesting
to note that  the additional factor $\ln^{\frac{d-1}{p}}(n)$ disappears in the dimension $d=1$. First, note that if $\alpha=0$ the one-dimensional setting
was considered in
\cite{Juditsky} and \cite{patricia}. The setting of \cite{Juditsky}
corresponds to $r=\infty$, while \cite{patricia} deal with the
case of $p=2$ and $\tau(2)>0$.
 Both settings rule out the sparse zone.
The rates of convergence found
 in these papers  are easily recovered from our results corresponding to the tail and dense zones.

 Next, we remark that the aforementioned factor appears only when anisotropic functional classes are considered.
Indeed, in view of the second assertion of Theorem \ref{th:adaptive-L_p-bounded-deconv}
our estimator is nearly optimally adaptive on the tail zone in the isotropic case. The natural question arising in this context, is whether the $\ln^{\frac{d-1}{p}}(n)$-factor is an unavoidable payment for anisotropy of the underlying function
or not?

At last, we note that in the isotropic case  our results remain true  when the corresponding Nikol'skii class is defined in $\bL_1$-norm
on $\bR^d$ ($\mathbf{r}=1$). It is worth noting that the analysis of the proof of the theorem allows us to assert that if $r_j=1$, $j\in J$ for some  $J\neq \{1,\ldots,d\}$ the first statement remains true up to some logarithmic factor. However the asymptotic of the maximal risk of our estimator if $r_j=1$ for any $j=1,\ldots,d$ remains unknown.

\vskip0.1cm

$\mathbf{5^0}.\;$ We finish our discussion with the following remark. If $\alpha\neq 1$ the assumption $f\in\bF_{g,\infty}(R,Q)$  implies in many cases that $f$ is uniformly bounded and, therefore, Theorem \ref{th:adaptive-L_p-bounded-deconv} is applicable. In particular it is always the case if the model (\ref{eq:observation-scheme}) is considered. Indeed $f,g\in\mP\big(\bR^d\big)$ in this case, which implies $\|f\|_\infty\leq (1-\alpha)^{-1}\|\mathfrak{p}\|_\infty\leq (1-\alpha)^{-1}Q$. Another case is  $\|g\|_\infty<\infty$ and recall that this assumption was used in the proofs of Theorems \ref{th1:lower-bound-deconvolution} and \ref{th2:lower-bound-deconvolution}, Assumption \ref{ass4:ass-on-noise-lower-bound}. We obviously have that
$$
\|f\|_\infty\leq (1-\alpha)^{-1}\big[Q+\alpha R\|g\|_\infty\big].
$$
More generally $\|f\|_\infty\leq (1-\alpha)^{-1}(Q+\alpha D)$ if $f\in\bF_{g,\infty}(R,Q)$ and $\|f\star g\|_\infty\leq D$. Since the definition of the Nikol'skii class implies that $\|f\|_{r^*}\leq L^*$, where $r^*=\sup_{j=1,\ldots, d}r_j$ and $L^*=\sup_{j=1,\ldots, d}L_j$, the latter condition can be verified in particular  if $\|g\|_q<\infty, 1/q=1-1/r^*$.
All saying above explains why we study the estimation of unbounded functions only in the case $\alpha=1$.

\subsubsection{\textsf{Unbounded case, $\alpha=1$}}
\label{sec:subsubsec:unbounded case}

 The   problem we address now  is the adaptive estimation over the collection of functional classes
$\big\{\bN_{\vec{r},d}\big(\vec{\beta},\vec{L}\big)\cap\bF_{g,\infty}(R,Q)\big\}_{\vec{\beta},\vec{r},\vec{L},R, Q}.\;$

\vskip0.05cm

As we already mentioned, if additionally $\|g\|_\infty<0$ then $\bF_{g,\infty}(R,Q)=\bF_{g}(R)$ for any $Q\geq R\|g\|_\infty$ and, therefore, in view of Theorem \ref{th1:lower-bound-deconvolution} discussed in Section \ref{sec:subsec-Lower bound for minimax-risk-deconv},
there is no consistent estimator if either $p=1$ or $\kappa_\alpha(p)\leq 0,\;\tau(p)\leq 0,\; p^*=p$. Analyzing the proof of the latter theorem, we come to the following assertion.

\begin{conjecture}
Let $\alpha=1$ and assume that Assumption \ref{ass3:ass-on-noise-lower-bound} is fulfilled. Suppose additionally that
Assumption \ref{ass1:ass-on-noise-lower-bound} holds with $\min_{j=1,\ldots,d} \mu_j>1/p$. Then, the assertion of Theorem \ref{th1:lower-bound-deconvolution} remains true if one replaces $\bN_{\vec{r},d}\big(\vec{\beta},\vec{L}\big)\cap\bF_{g}(R)$ by $\bN_{\vec{r},d}\big(\vec{\beta},\vec{L}\big)\cap\bF_{g,\infty}(R,Q)$.
\end{conjecture}

The latter result is formulated as a conjecture only because we will not prove it in the present paper. Its proof is postponed to Part III where the adaptive estimation over the collection
$$
\big\{\bN_{\vec{r},d}\big(\vec{\beta},\vec{L}\big)\cap\bF_{g,\mathbf{u}}(R,Q)\cap\bB_{\mathbf{q},d}(Q)\big\}_{\vec{\beta},\vec{r},\vec{L},R, Q, \mathbf{u},\mathbf{q}}
$$
introduced in Part I will be studied.
For this reason, later on we will only consider the parameters $\vec{\beta},\vec{r}$ belonging to the set $\cP_{p,\vec{\mu}}$ defined below.
\begin{equation*}
\label{eq:set-of-consistency-deconv}
\cP_{p,\vec{\mu}}=(0,\infty)^d\times[1,\infty]^d\setminus\Big\{\vec{\beta},\vec{r}:
\; \kappa_\alpha(p)\leq 0,\;\tau(p)\leq 0,\; \max_{j=1,\ldots,d}r_j\leq p\Big\}.
\end{equation*}
For given $p>1$ and $\vec{\mu}\in (0,\infty)^d$ the  latter set consists of the class parameters  for which a uniform consistent estimation is possible.
\begin{theorem}
\label{th:adaptive-L_p-general-deconv}
Let    $\ell\in\bN^*$ and $g\in\bL_1\big(\bR^d\big)$, satisfying Assumption \ref{ass1:ass-on-noise-upper-bound} be fixed
and let $K$ satisfy Assumptions \ref{ass:on-kernel-deconvolution} and \ref{ass2:on-kernel-deconvolution}.

\vskip0.1cm

1) Then for any $p>[\min_{j=1,\ldots} \mu_j]^{-1}$,  $R,Q>0$, $0<L_0\leq L_\infty<\infty$,  $\big(\vec{\beta},\vec{r})\in\cP_{p,\vec{\mu}}\cap \big\{(0,\ell]^d\times(1,\infty]^d\big\}$ and $\vec{L}\in [L_0,L_\infty]^d$
there exists $C<\infty$, independent of $\vec{L}$,  such that:
$$
\limsup_{n\to \infty}\sup_{f\in\bN_{\vec{r},d}\big(\vec{\beta},\vec{L}\big)\cap\bF_{g,\infty}(R,Q)}\mb_n\big(\cH^d\big)^{-1}\delta_n^{-\varrho(1)}
\cR^{(n)}_p\big[\widehat{f}_{\vec{\mathbf{h}},\cH^d}; f\big]\leq C,
$$
where  $\varrho(\cdot)$ is defined in (\ref{eq:varrho(alpha)}).

\vskip0.1cm

2) For any $p>[\min_{j=1,\ldots} \mu_j]^{-1}$, $R,Q>0$, $0<L_0\leq L_\infty<\infty$, $(\blb,\mathbf{r})\in \cP_{p,\vec{\mu}}\cap \big\{(0,\ell]\times(1,\infty]\big\}$ and $\mathbf{L}\in [L_0,L_\infty]$
there exists $C<\infty$, independent of $\mathbf{L}$,  such that:

$$
\limsup_{n\to \infty}\sup_{f\in\bN_{\mathbf{r},d}\big(\blb,\mathbf{L}\big)\cap\bF_{g,\infty}(R,Q)}\mb_n\big(\cH^d_{\text{isotr}}\big)^{-1}\delta_n^{-\varrho(1)}
\cR^{(n)}_p\big[\widehat{f}_{\vec{\mathbf{h}},\cH^d_{\text{isotr}}}; f\big]\leq C.
$$
\end{theorem}

Some remarks are in order.

$\mathbf{1^0}.\;$ Note that $\|g\|_1<\infty, \|g\|_\infty<\infty$ implies that $\|g\|_2<\infty$ and, therefore the Parseval identity together with  Assumption \ref{ass1:ass-on-noise-upper-bound} allows us to assert that
\begin{equation}
\label{eq:restriction-on-mu-part2}
\|g\|_\infty<\infty\quad \Rightarrow\quad \mu_j>1/2,\quad \forall j=1,\ldots,d.
\end{equation}
Hence,  the condition $p>[\min_{j=1,\ldots} \mu_j]^{-1}$ is automatically checked if $p\geq 2$ and $\|g\|_\infty<\infty$.

Also, it is worth noting that considering the adaptation over the collection of isotropic classes, we do not require that the coordinates of $\vec{\mu}$ would be the same. The latter is true for the second assertion of Theorem \ref{th:adaptive-L_p-bounded-deconv} as well. At last, analyzing the proof of the theorem, we can
assert that the second assertion remains true under the slightly weaker assumption $p>d(\mu_1+\cdots+\mu_d)^{-1}$.

\smallskip

$\mathbf{2^0}.\;$ The assertion of Theorem \ref{th1:lower-bound-deconvolution}   has no analogue in the existing literature  except the results obtained  in \cite{comte} and \cite{rebelles16}. \cite{comte} deals with the particular case  $p=2$, $\vec{r}=(2,\ldots,2)$ while \cite{rebelles16}
studied the case $\vec{r}=(p,\ldots,p)$, $p\in (1,\infty)$.
It is easy to check that in both papers whatever the value of $\vec{\beta}$ and $\vec{\mu}$, the corresponding set of parameters belongs to the \textsf{dense} zone. Note also that the estimation procedures used
in \cite{comte} as well as in \cite{rebelles16}, if $p\geq 2$, (both based on a global version of the Goldenshluger-Lepski method) are optimally-adaptive. They attain the asymptotic of  minimax risks corresponding to the dense zone found in Theorem \ref{th1:lower-bound-deconvolution}, while our method  is only nearly optimally adaptive.
However, it is well-known that the global selection from the family of standard kernel estimators
leads to correct results only if $\vec{r}=(p,\ldots,p)$ when the $\bL_p$-risk is considered, see, for instance \cite{GL11}.
On the other hand, estimation procedures based on a local selection scheme, which can be applied to the estimation of functions belonging to much more general functional classes, often do not lead to an optimally adaptive method.
Fortunately, the loss of accuracy inherent to local procedures is logarithmic w.r.t. the number of observations.

\smallskip

$\mathbf{3^0}.\;$ Together with Theorems \ref{th1:lower-bound-deconvolution} and  \ref{th2:lower-bound-deconvolution},  Theorems \ref{th:adaptive-L_p-bounded-deconv} and \ref{th:adaptive-L_p-general-deconv}   provide the  full classification of the asymptotics of the minimax risks over anisotropic/isotropic Nikolskii classes for the class parameters belonging to the sparse zone  and, up to some logarithmic factor, belonging to the tail and dense zones as well as the boundaries.
We mean that  the results of these theorems are valid for any \verb"fixed"
$\vec{\beta}\in (0,\infty)^d,\vec{r}\in(1,\infty]^d$ and $\vec{L}\in(0,\infty)^d$.  Indeed, for given $\vec{\beta}$ and $\vec{L}$ one can choose
$L_0=\min_{j=1,\ldots d}L_j$, $L_\infty=\max_{j=1,\ldots d}L_j$ and the number $\ell$, used in the  kernel construction (\ref{eq:w-function}), as any integer strictly larger than $\max_{j=1,\ldots d}\beta_j$.

\subsubsection{\textsf{Open problems}} Let us briefly discuss some unresolved adaptive estimation problems  in the convolution structure density model.


\paragraph{\textsf{Construction of an optimally-adaptive estimator}} As we already mentioned the proposed pointwise selection rule leads
to an optimal adaptive estimator only for the class parameters belonging to the sparse zone (in both bounded and unbounded case). We conjecture that
the construction of an optimally-adaptive estimator for all values of the nuisance parameters via pointwise selection is impossible, and other methods
should be invented. It is worth noting that no optimally-adaptive estimator is known  neither in the density model nor in the density
deconvolution even in dimension 1. In dimension larger than 1, one of the intriguing questions is related to the eventual price to pay for anisotropy
($\ln^{\frac{d-1}{p}}(n)$-factor) discussed in the remark $\mathbf{4^0}$ after Theorem \ref{th:adaptive-L_p-bounded-deconv}.

\paragraph{\textsf{Adaptive estimation of unbounded functions}} We were able to study the unbounded case only if $\alpha=1$. The estimation of unbounded densities under direct  as well as  partially contaminated observations remain open problems. We conjecture that the results obtained
in the case $\alpha=1$ are not true anymore for $\alpha\neq 1$ (neither upper bounds nor lower bound), but correct (or nearly correct) upper bounds for the asymptotics of the minimax risk can still be deduced from the oracle inequalities proved in Part I.

In the case  $\alpha=1$ there are at least two interesting problems. First,  all our results are valid under the condition $p>[\min_{j=1,\ldots} \mu_j]^{-1}$.  How the absence of this assumption may have effects on the accuracy of estimation is absolutely unclear.
Next,  let us mention that the lower bound result proved in Theorem \ref{th1:lower-bound-deconvolution} holds only under the consideration of the convolution
structure density model. Could the same bounds be established in the deconvolution model (\ref{eq:observation-scheme})?

\paragraph{\textsf{Adjustment  of "lower" and "upper bound" assumptions to each other}} Comparing the assertions of Theorems \ref{th1:lower-bound-deconvolution} and \ref{th2:lower-bound-deconvolution} with those of Theorem \ref{th:adaptive-L_p-bounded-deconv} and \ref{th:adaptive-L_p-general-deconv}, we remark that the obtention of the corresponding lower bounds for the minimax risk requires additional, rather restrictive, assumptions on the function $g$. Can they be weakened or even removed?

\section{Proof of Theorems \ref{th:adaptive-L_p-bounded-deconv} and \ref{th:adaptive-L_p-general-deconv}}
\label{sec:prof-of Theorem-th:adaptive-L_p-bounds-deconv}


The proofs are based on the application of Theorem 3 from Part I and on some auxiliary assertions presented below.

In the subsequent proof
$\mathbf{c},\mathbf{c}_1,\mathbf{c}_2, C, C_1, C_2 \ldots$,
stand for constants that can depend on $g, L_0,L_\infty$, $Q, R$,
$\vec{\beta}$, $\vec{r}$, $d$ and $p$, but
are independent of $\vec{L}$ and $n$. These constants can be different on
different appearances.

\subsection{Important concepts from Part I and proof outline}
\label{sec:subsec:Proof-outline}
In this section we recall the definition of some important quantities that appeared in Theorem 3 of Part I and discuss the facts which should be established to make this theorem applicable.

\vskip0.1cm

$\mathbf{I^0.\;}$  Theorem 3 (Part I)  deals with the minimax result over a class $\bF$ being an arbitrary subset of  $\bF_{g,\mathbf{u}} (R,D)\cap\bB_{\mathbf{q},d}(D)$ defined in Section 2.3 of Part I. In Theorem \ref{th:adaptive-L_p-bounded-deconv} we will consider $\bF=\bN_{\vec{r},d}\big(\vec{\beta},\vec{L}\big)\cap\bB_{\mathbf{\infty},d}(Q)$ and, therefore, $\bF\subset\bF_{g,\mathbf{\infty}}(R,D)\cap\bB_{\mathbf{\infty},d}(Q)$ with $D=Q[1-\alpha+\alpha\|g\|_1]$. This makes Theorem 3 (Part I) with $\mathbf{u}=\infty$  applicable in this case.

In Theorem \ref{th:adaptive-L_p-general-deconv} we consider
$\bF=\bN_{\vec{r},d}\big(\vec{\beta},\vec{L}\big)\cap\bF_{g,\mathbf{\infty}}(Q)$. We will show that for any $\vec{\beta},\vec{r}$ and $\vec{L}$ one can find
$\mathbf{q}>1$ and $D>0$ such that $\bN_{\vec{r},d}\big(\vec{\beta},\vec{L}\big)\subset\bB_{\mathbf{q},d}(D)$ and, therefore, Theorem 3 (Part I) is applicable with $\mathbf{u}=\infty$. The latter inclusions are mostly based on the embedding of anisotropic Nikol'skii spaces used in the proof of Proposition \ref{prop:measure-of-bias-deconv} and on Lemma \ref{lem:auxilary-results-for-prpositions-1-and2}.

\vskip0.1cm

$\mathbf{II^0.\;}$
The application of  Theorem 3 (Part I) in the case $\mathbf{u}=\infty$ requires to compute
\begin{gather*}
J\big(\vec{h},v\big)=\big\{j\in\{1,\ldots,d\}:\; h_j\in\mathbf{V}_{j}(v)\big\},\;\;\mathbf{V}_{j}(v)=\big\{\mathbf{v}\in\cH:\;\; \mathbf{B}_{j,\infty,\bF}(\mathbf{v})\leq \mathbf{c}v\big\},
\nonumber\\
\label{eq4:def-blL}
\blL_{\vec{s}}(v,\bF,\mathbf{\infty})=\inf_{\vec{h}\in\mH(v,2)}\bigg[\sum_{j\in\bar{J}(\vec{h},v)} v^{-s_j}\big[\mathbf{B}_{j,s_j,\bF}\big(h_{j}\big)\big]^{s_j}\bigg];
\\
\label{eq1:def-blL}
\blL_{\vec{s}}\big(v,\bF\big)=\inf_{\vec{h}\in\mH(v)}\bigg[\sum_{j\in\bar{J}(\vec{h},v)} v^{-s_j}\big[\mathbf{B}_{j,s_j,\bF}\big(h_{j}\big)\big]^{s_j}+v^{-2}F_n^2\big(\vec{h}\big)\bigg],
\end{gather*}
where remind
$
F_n\big(\vec{h}\big)=\big(\ln{n}+ \sum_{j=1}^d|\ln{h_j|}\big)^{1/2}\prod_{j=1}^d(nh)^{-\frac{1}{2}}_j(h_j\wedge 1)^{-\blg_j(\alpha)}
$
and $\mathbf{c}>0$ is a universal constant completely determined by the kernel $\cK_\ell$ and the dimension $d$.

 In the next section we propose  quite sophisticated constructions of vectors $\blh(\cdot,\mathbf{s})$ and
 $\vec{\mh}(\cdot,\mathbf{s})$, $\mathbf{s}\in[1,\infty]$ and show, Propositions \ref{prop1:auxiliary results-new} and \ref{prop2:auxiliary results-new}, that
\begin{equation}
\label{eq1:preliminaries-proof-theorems}
 \vec{\blh}(v,\mathbf{1})\in \mH(v),\; v\in[\underline{\mathbf{v}},1],\quad \vec{\blh}(v,\mathbf{\infty})\in \mH(v,2),\;v\in[\Bv,\overline{\Bv}],\quad \vec{\mh}(\mathbf{v},\mathbf{\infty})\in \mH(\mathbf{v},2).
 \end{equation}
 Here  $\Bv$ is defined in (\ref{eq02:new-paper}), $\underline{\mathbf{v}},\mathbf{v}$ are defined in (\ref{eq0200:new-paper}) and  $\overline{\Bv}\in \{1,\mathbf{v_1},\mathbf{v_3},\overline{\mathbf{v}}, \overline{\mathbf{v}}\wedge\mathbf{v_3}\}$, where $\mathbf{v_1}, \mathbf{v_3}$ are defined in (\ref{eq0201:new-paper}) and $\overline{\mathbf{v}}$ is given in (\ref{eq02010:new-paper}).
In  Proposition \ref{prop:measure-of-bias-deconv} we prove that for any $\vec{h}\in\cH^d$
\begin{equation}
\label{eq2:preliminaries-proof-theorems}
\mathbf{B}_{j,r_j,\bN_{\vec{r},d}\big(\vec{\beta},\vec{L}\big)}\big(h_{j}\big)\leq C_1L_jh_j^{\beta_j},\quad j=1,\ldots,d.
\end{equation}
and if $\tau(p^*)>0$ then  additionally
 \begin{equation}
\label{eq3:preliminaries-proof-theorems}
\mathbf{B}_{j,q_j,\bN_{\vec{r},d}\big(\vec{\beta},\vec{L}\big)}\big(h_{j}\big)\leq C_1L_jh_j^{\gamma_j},\quad j=1,\ldots,d,
\end{equation}
where $\vec{\gamma}$ and $\vec{q}$ are defined in (\ref{eq:gamma-and-q}) below and $C_1$ is independent of $\vec{L}$.
At last the definition of $\blh(\cdot,\mathbf{s})$ and
 $\vec{\mh}(\cdot,\mathbf{s})$, $\mathbf{s}\in[1,\infty]$ together with (\ref{eq2:preliminaries-proof-theorems}) allows us to assert, see  (\ref{eq91:preliminaries-proof-theorems}), that
 \begin{equation}
\label{eq4:preliminaries-proof-theorems}
J\big(\vec{\blh}(v,\mathbf{1}),v\big)\supseteq \cJ_\infty,\;\;\;
J\big(\vec{\blh}(v,\mathbf{\infty}),v\big)\supseteq \cJ_\infty,\;\;\;
J\big(\vec{\mh}(v,\mathbf{\infty}),v\big)\supseteq \cJ_\infty,\quad \forall v>0,
\end{equation}
where $\cJ_\infty=\{j=1,\ldots,d:\; r_j=\infty\}$. Thus, putting
\begin{gather*}
\label{eq5:preliminaries-proof-theorems}
\bll_{1}(v)=\sum_{j\in\bar{\cJ}_\infty} v^{-r_j}L^{r_j}_j\big[\blh_j(v,\mathbf{\infty})\big]^{r_j\beta_j},\quad
\bll=\sum_{j\in\bar{\cJ}_\infty} \mathbf{v}^{-q_j}L^{q_j}_j\big[\mh_j(\mathbf{v},\mathbf{\infty})\big]^{q_j\gamma_j};
\nonumber\\
\label{eq6:preliminaries-proof-theorems}
\bll_{2}(v)=\sum_{j\in\bar{\cJ}_\infty} v^{-r_j}L^{r_j}_j\big[\blh_j(v,\mathbf{1})\big]^{r_j\beta_j}+v^{-2}(\ln{n}/n)
\prod_{j=1}^d(\blh_j(v,\mathbf{1}))^{-1-2\blg_j(\alpha)},
\end{gather*}
 we obtain in view of  (\ref{eq1:preliminaries-proof-theorems}), (\ref{eq2:preliminaries-proof-theorems}) and (\ref{eq4:preliminaries-proof-theorems}) that
\begin{eqnarray}
\label{eq7:preliminaries-proof-theorems}
\blL_{\vec{r}}\big(v,\bN_{\vec{r},d}\big(\vec{\beta},\vec{L}\big),\mathbf{\infty}\big)&\leq& C_1\bll_{1}(v),\quad \forall v\in[\Bv,\overline{\Bv}];
\\
\label{eq8:preliminaries-proof-theorems}
 \blL_{\vec{r}}\big(v,\bN_{\vec{r},d}\big(\vec{\beta},\vec{L}\big)\big)&\leq& C_1\bll_{2}(v), \quad \forall v\in[\underline{\mathbf{v}},1]
\end{eqnarray}
To get (\ref{eq8:preliminaries-proof-theorems}) we have used that for all $n$ large enough and all $ v\in[\underline{\mathbf{v}},1]$

\vskip0.2cm

\centerline{$
F_n\big(\vec{\blh}(v,\mathbf{1}))\big)\leq C_2(\ln{n}/n)
\prod_{j=1}^d(\blh_j(v,\mathbf{1}))^{-1-2\blg_j(\alpha)},
$}

\vskip0.2cm

\noindent where $C_2$ is independent of $\vec{L}$. This follows from assertions (\ref{eq01:proof-prop1:auxiliary results-new}) and (\ref{eq1:proof-prop1:auxiliary results-new}) established in the proof of Proposition \ref{prop1:auxiliary results-new}.
We deduce from  (\ref{eq7:preliminaries-proof-theorems}) and  (\ref{eq8:preliminaries-proof-theorems}),
 the following bound.
\begin{gather}
\label{eq10:preliminaries-proof-theorems}
\big[\blL_{\vec{r}}\big(v,\bN_{\vec{r},d}\big(\vec{\beta},\vec{L}\big),\mathbf{\infty}\big)
\wedge\blL_{\vec{r}}\big(v,\bN_{\vec{r},d}\big(\vec{\beta},\vec{L}\big)\big)\big]\leq C_1\Big[\bll_2(v)\mathrm{1}_{[\underline{\Bv},\Bv]}(v)+ \bll_1(v)\mathrm{1}_{[\Bv,\overline{\Bv}]}(v)\Big].
\end{gather}
Moreover, if $\tau(p^*)>0$ we get in view of (\ref{eq1:preliminaries-proof-theorems}), (\ref{eq3:preliminaries-proof-theorems}) and (\ref{eq4:preliminaries-proof-theorems})
\begin{gather}
\label{eq9:preliminaries-proof-theorems}
\mathbf{v}^p\blL_{\vec{q}}\big(\mathbf{v},\bN_{\vec{r},d}\big(\vec{\beta},\vec{L}\big),\mathbf{\infty}\big)\leq C_1\mathbf{v}^p\bll.
\end{gather}

\subsection{Special set of bandwidths}
\label{sec:subsec-band-construction}
The bandwidth's construction presented below as well as auxiliary statements from the next section will be exploited not only for proving Theorems \ref{th:adaptive-L_p-bounded-deconv} and \ref{th:adaptive-L_p-general-deconv}, but also in the consideration forming Part III of this work. By this reason we formulate them in a bit more general form than what is needed for our current purposes.
Set for any $r,\mathbf{s}\in [1,\infty]$
\vskip0.15cm
\centerline{$
\kappa_\alpha(r,\mathbf{s})=\frac{\mathbf{s}\omega(\alpha)(2+1/\beta(\alpha))}{(\mathbf{s}+\omega(\alpha))}-r,\quad \alpha\in[0,1].
$}

\vskip0.1cm

\noindent Recall that $\mathbf{c}= \big(20d\big)^{-1}\big[\max (2c_{\cK_\ell}\|\cK_\ell\|_\infty,\|\cK_\ell\|_1)\big]^{-d}$ and let $\BL>0$ be any number satisfying (recall that $C_1$ appeared in (\ref{eq2:preliminaries-proof-theorems}))
\begin{equation}
\label{eq90:preliminaries-proof-theorems}
\BL\leq 1\wedge (C_1^{-1}\mathbf{c}) \wedge L_0.
\end{equation}
Recall that
$
\delta_n=L(\alpha)n^{-1}\ln{n}
$
and
introduce for any   $v>0$, $\mathbf{s}\in [1,\infty]$ and   $j=1,\ldots,d$
\begin{eqnarray}
\label{eq3::def-bar(mu)_j}
\widetilde{\bleta}_j(v,\mathbf{s})&=&
\big(\BL L_j^{-1}\big)^{\frac{1}{\beta_j}}\big\{\ma^{-2}\delta_n\big\}^{\frac{\mathbf{s}\omega(\alpha)}{(\mathbf{s}+\omega(\alpha))\beta_jr_j}}
 v^{\frac{1}{\beta_j}-\frac{\mathbf{s}\omega(\alpha)(2+1/\beta(\alpha))}{(\mathbf{s}+\omega(\alpha))\beta_jr_j}};
 \\
\widehat{\bleta}_j(v,\mathbf{s})&=&\big(\BL L_j^{-1}\big)^{\frac{1}{\gamma_j}}
\big\{\ma^{-2}\delta_n\big\}^{\frac{\mathbf{s}\upsilon(\alpha)}{(\mathbf{s}+\upsilon(\alpha))\gamma_jq_j}}
 v^{\frac{1}{\gamma_j}-\frac{\mathbf{s}\upsilon(\alpha)(2+1/\gamma(\alpha))}{(\mathbf{s}+\upsilon(\alpha))\gamma_jq_j}},
\end{eqnarray}
where we have put $p_\pm=[\sup_{j\in\bar{\cJ}_\infty}r_j]\vee p$, $\bar{\cJ}_\infty$ is complimentary to $\cJ_\infty$  and
\begin{gather}
\label{eq:gamma-and-q}
\qquad\;\; q_j=\left\{\begin{array}{ll}
p_\pm,\quad & j\in\bar{\cJ}_\infty,
\\
\infty,\quad & j\in\cJ_\infty,
\end{array}
\right.,
\qquad\;\;\;\;
\gamma_j=\left\{
\begin{array}{ll}
\frac{\beta_j\tau(p_\pm)}{\tau(r_j)},\quad & j\in\bar{\cJ}_\infty,
\\
\beta_j,\quad & j\in \cJ_\infty.
\end{array}
\right.
\\
\nonumber
\label{eq:gamma-upsilon}
\frac{1}{\gamma(\alpha)}:=\sum_{j=1}^d \frac{2\blg_j(\alpha)+1}{\gamma_j},\qquad\quad
\frac{1}{\upsilon(\alpha)}:=\sum_{j=1}^d \frac{2\blg_j(\alpha)+1}{\gamma_j q_j}.
\end{gather}
The constant $\ma>0$  will be chosen differently in accordance with  some special relationships between the parameters $\vec{\beta}$, $\vec{r}$, $\vec{\mu}$, $\alpha$ and $p$.
Determine $\blh_j(\cdot,\mathbf{s})$ and $\mh_j(\cdot,\mathbf{s}), j=1,\ldots,d$, from the relations
\begin{eqnarray}
\label{eq:def-h(v)}
\blh_j(v,\mathbf{s})&=&\max\big\{h\in\cH:\;\; h\leq \widetilde{\bleta}_j(v,\mathbf{s})\big\},\;\;v>0;
\\*[2mm]
\mh_j(v,\mathbf{s})&=&\max\big\{h\in\cH:\;\; h\leq \widehat{\bleta}_j(v,\mathbf{s})\big\},\;\;v>0,
\end{eqnarray}
and set $\vec{\blh}(\cdot,\mathbf{s})=\big(\blh_1(\cdot,\mathbf{s}),\ldots, \blh_d(\cdot,\mathbf{s})\big)$ and $\vec{\mh}(\cdot,\mathbf{s})=\big(\mh_1(\cdot,\mathbf{s}),\ldots,\mh_d(\cdot,\mathbf{s})\big)$.

\subsection{Auxiliary statements} All the results formulated below are proved in Section \ref{sec:Proofs-of-Propositions}. Let

\vskip0.2cm

\centerline{$
\mz(v)=2\big(\ma^{-2}\delta_n\big)^{-\frac{\omega(\alpha)}{\omega(\alpha)+\mathbf{u}}}
v^{\frac{\omega(\alpha)(2+1/\beta(\alpha))}{\mathbf{u}+\omega(\alpha)}},\quad \mathbf{u}\in[1,\infty],
$}

\vskip0.2cm

\noindent and remark that  $\mz(\cdot)\equiv2$ if $\mathbf{u}=\infty$.
Note also that
\begin{equation}
\label{eq02:new-paper}
\mz(v)\geq 2,\quad \forall v\geq \big(\ma^{-2}\delta_n\big)^{\frac{1}{2+1/\beta(\alpha)}}=\Bv.
\end{equation}
Introduce the following notations: $\mu(\alpha)=\min_{j=1,\ldots,d}\mu_j(\alpha)$,
$$
X=\frac{1}{2\beta(1)}-\frac{1}{2\beta(0)}=\sum_{j=1}^d\frac{\blg_j(\alpha)}{\beta_j},\;\;\;
Y=\frac{1}{2\omega(1)}-\frac{1}{2\omega(0)}=\sum_{j=1}^d\frac{\blg_j(\alpha)}{\beta_jr_j}.
$$
Recall that  $z(\alpha)=\omega(\alpha)
(2+1/\beta(\alpha))\beta(0)\tau(\infty)+1$ and define
\begin{equation}
\label{eq0200:new-paper}
\underline{\mathbf{v}}=(\ma^{-2}\delta_n)^{\frac{1}{1-1/\omega(\alpha)+1/\beta(\alpha)}}, \quad \mathbf{v}=\big(\ma^{-2}\delta_n\big)^{\frac{\omega(\alpha)\tau(\infty)\beta(0)}{z(\alpha)+\omega(\alpha)/\mathbf{u}}}.
\end{equation}
Set  $\mathbf{u}^*=[-\tau(\infty)\beta(0)]^{-1}$ if $\tau(\infty)<0$ and let $\mathbf{u}^*=\infty$ if $\tau(\infty)\geq 0$. Put finally $\mathbf{y}=\mathbf{u}^*\vee p^*$.

\begin{proposition}
\label{prop1:auxiliary results-new}
Let $\vec{\beta}$, $\vec{r}$, $L_0,L_\infty$, $\vec{\mu}$, $\alpha$ and $p$ be given. Assume that $\vec{L}\in[L_0,L_\infty]^d$. Then,

1) there exists $\ma>0$ independent of $\vec{L}$ such that for all $n$ large enough
$$
\vec{\blh}(v,\mathbf{1})\in \mH(v),\quad \forall v\in [\underline{\mathbf{v}},1],
$$
\quad 2) there exists $\ma>0$ independent of $\vec{L}$ and $\mathbf{u}$ such that for all $n$ large enough
$$
\vec{\mh}(\mathbf{v},\mathbf{u})\in \mH\big(\mathbf{v},\mz(\mathbf{v})\big)
$$
if either $\kappa_\alpha(p^*,\mathbf{u})<0, \tau(\infty)\geq 0$ or  $\kappa_\alpha(p^*,\mathbf{u})<0,$ $\tau(p^*)>0$, $Y\geq[X+1]\mathbf{y}^{-1}-1/\mathbf{u}$.

\end{proposition}

\begin{remark}
\label{rem:after-prop1}
Note that if $\alpha\neq1$, the condition $Y \ge [X+1]\mathbf{y}^{-1}-1/\mathbf{u}$ simply means $\mathbf{u}\leq \mathbf{u^*}\vee p^*$, since $X=Y=0$.
On the other hand if $\alpha=1$ this condition holds if $\tau(\infty)\geq 0$ whatever the values of $\vec{\beta},\vec{\mu}$ and $\vec{r}$, since $Y>0$.
Also, note that
\begin{equation}
\label{eq:rem:about-main-condition}
\mu(1)+1/\mathbf{u}-1/\mathbf{y}\geq0\quad\Rightarrow\quad Y\geq[X+1]\mathbf{y}^{-1}-1/\mathbf{u}.
\end{equation}
Indeed, since $r_j\leq p^*\leq \mathbf{y}$ for any $j=1,\ldots,d$ we have
$$
 Y-[X+1]\mathbf{y}^{-1}+1/\mathbf{u}\geq \mu(1)[1-\tau(\mathbf{y})]-1/\mathbf{y}+1/\mathbf{u}\geq  \mu(1)-1/\mathbf{y}+1/\mathbf{u}
$$
and (\ref{eq:rem:about-main-condition}) follows. To get the last inequality we have used that $\tau(\mathbf{u}^*)=0$ and that $\tau(\cdot)$ is strictly decreasing, so $\tau(\mathbf{y})\le 0$.
In particular we deduce from (\ref{eq:rem:about-main-condition}) that the condition $ Y>[X+1]\mathbf{y}^{-1}-1/\mathbf{u}$  is always fulfilled in the case $\mathbf{u}=\mathbf{u}^*$.

\end{remark}

Recall that $\Bv\to 0,n\to\infty,$ is defined in (\ref{eq02:new-paper}) and introduce the following quantities.
\begin{gather}
\label{eq0201:new-paper}
\mathbf{v_1}=\big(\ma^{-2}\delta_n\big)^{\frac{1}{1-\mathbf{u}/\omega(0)+1/\beta(0)}},\;\; \mathbf{v_2}=\big(\ma^{-2}\delta_n\big)^{\frac{\mathbf{u}\omega(1)}{\kappa_1(p^*,\mathbf{u})(\omega(1)+\mathbf{u})}},\;\;
\mathbf{v_3}=
\big(\ma^{-2}\delta_n\big)^{-\frac{Y+1/\mathbf{u}}{\pi(\mathbf{u})\vee 0}},
\end{gather}
where $\pi(\mathbf{u)}=[1/\omega(0)-1/\mathbf{u}][1+X]-1/\beta(0)[Y+1/\mathbf{u}].\;$  Define  also
\begin{equation}
\label{eq02010:new-paper}
\overline{\mathbf{v}}=\mathbf{v}\mathrm{1}_{\{\tau(p^*)>0\}}+\mathbf{v_2}\mathrm{1}_{\{\tau(p^*)\leq0\}}
\end{equation}
Note that $\mathbf{v_1}\to\infty, n\to\infty$, if
$\infty>\mathbf{u}\geq \mathbf{u}^*\vee p^*$ (it will be proved in Proposition \ref{prop2:auxiliary results-new} below). However $\mathbf{v_1}=1$ if $\mathbf{u}=\infty$. As it is shown in the proof of Proposition \ref{prop1:auxiliary results-new}, formulae (\ref{eq503:proof-prop1:auxiliary results-new}), $\Bv<\mathbf{v}$ for all $n$  large enough.
Also $\mathbf{v_2}\to\infty, n\to\infty$, if $\kappa_1(p^*,\mathbf{u})< 0$.
At last $\mathbf{v_3}\to\infty, n\to\infty$,  since $\omega(0)>\omega(1)$. Moreover $\mathbf{v_3}=\infty$ if $\pi(\mathbf{u})\leq 0$.
Introduce finally
\begin{eqnarray*}
\cI_{\mathbf{u}}(\alpha)&=&\left\{
\begin{array}{lllll}
\;[\Bv,1], &p^*=\infty
\\*[2mm]
\;[\Bv,\mathbf{v}_1],\; &\alpha\neq 1,\; p^*<\infty;
\\*[2mm]
\; [\Bv,\mathbf{v}_3], \; &\alpha=1,\; p^*<\infty,\;\kappa_\alpha(p^*,\mathbf{u})\geq 0;
\\*[2mm]
\;[\Bv,\overline{\mathbf{v}}], \; &\alpha=1,\; p^*<\infty,\;\kappa_\alpha(p^*,\mathbf{u})< 0,\; Y\geq[X+1]\mathbf{y}^{-1}-1/\mathbf{u};
\\*[2mm]
\;[\Bv,\overline{\mathbf{v}}\wedge\mathbf{v_3}], \; &\alpha=1,\; p^*<\infty,\;\kappa_\alpha(p^*,\mathbf{u})< 0,\; Y<[X+1]\mathbf{y}^{-1}-1/\mathbf{u},
\end{array}
\right.
\end{eqnarray*}

\begin{proposition}
\label{prop2:auxiliary results-new} Let $\vec{\beta}$, $\vec{r}$, $L_0,L_\infty$, $\vec{\mu}$, $\alpha$ and $p$ be given and let $\vec{L}\in[L_0,L_\infty]^d$, $\mathbf{u}\in [\mathbf{u}^*\vee p^*,\infty]$.
Then, there exists $\ma>0$  independent of $\vec{L}$ and $\mathbf{u}$  such that for all $n$ large enough
$$
\vec{\blh}(v,\mathbf{u})\in \mH\big(v,\mz(v)\big), \quad v\in\cI_{\mathbf{u}}(\alpha).
$$

\end{proposition}
In the current paper we will use the statements of Proposition \ref{prop1:auxiliary results-new} and \ref{prop2:auxiliary results-new} only with $\mathbf{u}=\infty$. In this context we remark that  $\kappa_\alpha(\cdot)\equiv\kappa_\alpha(\cdot,\mathbf{\infty})$.

\begin{proposition}
\label{prop:measure-of-bias-deconv}
Let $\ell\in\bN^*$, $p>1$ and  $K$ satisfying Assumption \ref{ass2:on-kernel-deconvolution} be fixed. Then for any  $\vec{\beta}\in(0,\ell]^d$, $\vec{r}\in [1,\infty]^d$ and $\vec{L}\in (0,\infty)^d$
 one can find $C_1>0$ independent of $\vec{L}$ such that
(\ref{eq2:preliminaries-proof-theorems}) holds. If additionally
 $\tau(p^*)>0$ then (\ref{eq3:preliminaries-proof-theorems}) is fulfilled as well. At last, (\ref{eq2:preliminaries-proof-theorems}) and (\ref{eq3:preliminaries-proof-theorems}) remain true if one replaces the quantity $\mathbf{B}$ by $\mathbf{B}^*$.

\end{proposition}

The quantities $\mathbf{B}_{j,s,\bF}(\cdot)$ and $\mathbf{B}^*_{j,s,\bF}(\cdot)$ are introduced in Part I but the reader can find them in the proof of the proposition.
Let us also present the following auxiliary results which will be useful in the sequel. Their proofs are postponed to Appendix.

\begin{lemma}
\label{lem:auxilary-results-for-prpositions-1-and2}
For any $\mathbf{u}\in[1,\infty]$
\begin{eqnarray}
\label{eq4:proof-prop1:auxiliary results-new}
\kappa_\alpha(p^*,\mathbf{u})\leq 0,\; \tau(p^*)>0,\quad&\Rightarrow&\quad z(\alpha)+\omega(\alpha)/\mathbf{u}>0;
\\
\label{eq900:proof-prop1:auxiliary results-new}
\quad Y\geq[X+1]\mathbf{y}^{-1}-1/\mathbf{u},\; \tau(p^*)>0,\quad&\Rightarrow&\quad z(\alpha)/\omega(\alpha)-1+2/\mathbf{u}\geq 0.
\end{eqnarray}
Let $Y-[X+1]\mathbf{y}^{-1}> 0$ and $\kappa_1(p^*,\mathbf{\infty})\geq 0$. Then there exists $s>p^*$ such that
\begin{equation}
\label{eq3:new-formulas}
\tau(s)>0,\quad s\geq (1+X)/Y.
\end{equation}

\end{lemma}

\noindent We finish this section with the following observations which will be useful in the sequel.

If $\kappa_\alpha(p^*)\geq 0$ one has
\begin{equation}
 \label{eq-new-new:nu}
 \varrho(\alpha)=\frac{1-1/p}{1-1/\omega(\alpha)+1/\beta(\alpha)}\bigwedge\frac{\beta(\alpha)}{2\beta(\alpha)+1}:=r(\alpha),\; \quad
 \rho(\alpha)=r(\alpha)\bigwedge \frac{\omega(\alpha)}{p} .
\end{equation}
\quad If $\kappa_\alpha(p^*)<0$ one has
\begin{eqnarray}
\label{eq-new:nu}
\varrho(\alpha)&=&r(\alpha)\bigwedge\bigg[\frac{\tau(p)\omega(\alpha)\beta(0)}{z(\alpha)}\mathrm{1}_{\{\tau(p^*)> 0\}}+
\frac{\omega(\alpha)(1-p^*/p)}{\kappa_\alpha(p^*)}\mathrm{1}_{\{\tau(p^*)\leq 0\}}\bigg];
\\
\label{eq:nu}
\rho(\alpha)&=&r(\alpha)\bigwedge\bigg[\frac{\tau(p)\omega(\alpha)\beta(0)}{z(\alpha)}\mathrm{1}_{\{\tau(\infty)> 0\}}+
\frac{\omega(\alpha)}{p}\mathrm{1}_{\{\tau(\infty)\leq 0\}}\bigg].
\end{eqnarray}

\subsection{Concluding remarks}
\label{sec:subsec-conluding-remarks}
Let us collect some bounds for several terms appearing in Theorem 3 (Part I) and used in the proofs of Theorems \ref{th:adaptive-L_p-bounded-deconv} and \ref{th:adaptive-L_p-general-deconv} simultaneously.

$\mathbf{1^0.\;}$First we remark that $\blh_j(\cdot,\mathbf{1})\equiv \blh_j(\cdot,\mathbf{\infty})\equiv\mh_j(\cdot,\mathbf{\infty})\leq \big(\BL L_j^{-1}\big)^{\frac{1}{\beta_j}}$, $j\in\cJ_\infty$.
Then, (\ref{eq4:preliminaries-proof-theorems}) follows from   (\ref{eq2:preliminaries-proof-theorems}) and (\ref{eq90:preliminaries-proof-theorems}) because for any $j\in\cJ_\infty$ and $v>0$
\begin{eqnarray}
\label{eq91:preliminaries-proof-theorems}
&&\;\;\;\mathbf{B}_{j,\infty,\bN_{\vec{r},d}\big(\vec{\beta},\vec{L}\big)}\big(\blh_{j}(v,\mathbf{1})\big)=
\mathbf{B}_{j,\infty,\bN_{\vec{r},d}\big(\vec{\beta},\vec{L}\big)}\big(\blh_{j}(v,\mathbf{\infty})\big)
=\mathbf{B}_{j,\infty,\bN_{\vec{r},d}\big(\vec{\beta},\vec{L}\big)}\big(\mh_{j}(v,\mathbf{\infty})\big)
\leq \mathbf{c}v.
\end{eqnarray}

$\mathbf{2^0.\;}$ We deduce from the definition of $\vec{\blh}(\cdot,\mathbf{s}), \mathbf{s}\in\{1,\infty\}$ that
$$
\bll_1(v)\leq c_1 \delta_n^{\omega(\alpha)}
 v^{-\omega(\alpha)(2+1/\beta(\alpha))},\;v\in\cI_\infty(\alpha), \quad \bll_2(v)\leq c_1 \delta_n^{\frac{\omega(\alpha)}{\omega(\alpha)+1}}
 v^{-\frac{\omega(\alpha)(2+1/\beta(\alpha))}{\omega(\alpha)+1}},\;v\in[\underline{\mathbf{v}},1].
$$
It yields together with (\ref{eq10:preliminaries-proof-theorems}) and the definitions of $\underline{\mathbf{v}}$ and $\Bv$, choosing $\underline{\Bv}=\underline{\mathbf{v}}$,
\begin{eqnarray}
\label{eq12:preliminaries-proof-theorems}
&&\hskip-0.8cm\int_{\underline{\Bv}}^{\overline{\Bv}} v^{p-1} \big[\blL_{\vec{r}}\big(v,\bN_{\vec{r},d}\big(\vec{\beta},\vec{L}\big),\mathbf{\infty}\big)
\wedge\blL_{\vec{r}}\big(v,\bN_{\vec{r},d}\big(\vec{\beta},\vec{L}\big)\big)\big]\rd v \leq c_2 \bigg[ \delta_n^{\frac{\omega(\alpha)}{\omega(\alpha)+1}}
 \underline{\mathbf{v}}^{p-\frac{\omega(\alpha)(2+1/\beta(\alpha))}{\omega(\alpha)+1}}\mathrm{1}_{\{\kappa_\alpha(p)>p\omega(\alpha)\}}
\nonumber\\
&&\qquad+\delta_n^{\frac{\omega(\alpha)}{\omega(\alpha)+1}}
 \Bv^{p-\frac{\omega(\alpha)(2+1/\beta(\alpha))}{\omega(\alpha)+1}}
 \mathrm{1}_{\{\kappa_\alpha(p)<p\omega(\alpha)\}}+\delta_n^{\omega(\alpha)}
 \Bv^{p-\omega(\alpha)(2+1/\beta(\alpha))}\mathrm{1}_{\{\kappa_\alpha(p)>0\}}
 \nonumber\\
 &&\qquad +\delta_n^{\omega(\alpha)}
 \overline{\Bv}^{p-\omega(\alpha)(2+1/\beta(\alpha))}\mathrm{1}_{\{\kappa_\alpha(p)<0\}}
 +\ln{(n)}\Big(\delta_n^{\frac{\omega(\alpha)}{\omega(\alpha)+1}}\mathrm{1}_{\{\kappa_\alpha(p)=p\omega(\alpha)\}}+
 \delta_n^{\omega(\alpha)}\mathrm{1}_{\{\kappa_\alpha(p)=0\}}\Big)\bigg]
 \nonumber\\
&&\qquad =: A_n+c_2 \delta_n^{\omega(\alpha)}\overline{\Bv}^{p-\omega(\alpha)(2+1/\beta(\alpha))}\mathrm{1}_{\{\kappa_\alpha(p)<0\}}.
\end{eqnarray}
After elementary computations and taking into account (\ref{eq-new-new:nu}), we obtain
\begin{eqnarray}
\label{eq13:preliminaries-proof-theorems}
A_n\leq c_3\mb^p_n(\bH)\delta_n^{p\rho(\alpha)},\qquad A_n\leq c_3\mb^p_n(\bH)\delta_n^{p\varrho(\alpha)}.
\end{eqnarray}
These bounds are not surprising because $\varrho(\alpha)=\rho(\alpha)$ if $\kappa_\alpha(p)\geq 0$.
At last, if $\tau(p^*)>0$, we get from (\ref{eq9:preliminaries-proof-theorems}) thanks to the definition of $\vec{\mh}(\cdot,\mathbf{\infty})$ and the presentation proved in (\ref{eq50:proof-prop1:auxiliary results-new})
with $\mathbf{u}=\infty$
\begin{gather}
\label{eq95:preliminaries-proof-theorems}
\mathbf{v}^p\blL_{\vec{q}}\big(\mathbf{v},\bN_{\vec{r},d}\big(\vec{\beta},\vec{L}\big),\mathbf{\infty}\big) \leq c_4
 \delta_n^{\frac{\omega(\alpha)\tau(p)\beta(0)}{z(\alpha)}}.
\end{gather}

$\mathbf{3^0.\;}$ At last, choosing $\underline{\Bv}=\underline{\mathbf{v}}$, we obtain $\ell_{\bH}(\underline{\Bv})\leq c_6 \delta_n^{\frac{p-1}{1-1/\omega(\alpha)+1/\beta(\alpha)}}\big(\ln{n}\big)^{t(\bH)}$, which yields by (\ref{eq-new-new:nu}), (\ref{eq-new:nu}) and (\ref{eq:nu}):
\begin{gather}
\label{eq100:preliminaries-proof-theorems}
\ell_{\bH}(\underline{\Bv})\leq c_6 \mb^p_n(\bH)\delta_n^{p \rho (\alpha)},\qquad
\ell_{\bH}(\underline{\Bv})\leq c_6 \mb^p_n(\bH)\delta_n^{p \varrho (\alpha)}.
\end{gather}

\subsection{Proof of Theorem \ref{th:adaptive-L_p-bounded-deconv}}
As it has already been mentioned we will apply Theorem 3 (Part I) with $\mathbf{u}=\infty$, $\mathbf{q}=\infty$, $D=Q[1-\alpha+\alpha\|g\|_1]\vee Q$
and $\underline{\Bv}=\underline{\mathbf{v}}$.

\vskip0.1cm

$\mathbf{1^0.\;}$ Consider the cases $\kappa_\alpha(p^*)\geq 0, $ or $\kappa_\alpha(p^*)<0, \tau(\infty)\leq0$.

Choose $\overline{\Bv}=1$ and remark that the statements of Propositions \ref{prop1:auxiliary results-new} and \ref{prop2:auxiliary results-new} hold for any
$v\in[\underline{\Bv},\overline{\Bv}]$.
Indeed, it suffices to note that $\cI_{\mathbf{\infty}}(\alpha)\supseteq[\underline{\Bv},\overline{\Bv}]:=[\underline{\mathbf{v}},1]$, because
$\mathbf{v_1},\mathbf{v_2},\mathbf{v_3}>1$ and  $\overline{\mathbf{v}}\geq 1$ if $\tau(\infty)<0$ since in this case $\mathbf{v}>1$ by (\ref{eq4:proof-prop1:auxiliary results-new}).  Then we can apply all the bounds obtained above, and in particular we get from (\ref{eq7:preliminaries-proof-theorems})
\begin{gather}
\label{eq96:preliminaries-proof-theorems}
 \blL_{\vec{r}}\big(1,\bN_{\vec{r},d}\big(\vec{\beta},\vec{L}\big),\mathbf{\infty}\big)\leq C_1\bll_1(1)\leq c_5\delta_n^{\omega(\alpha)}\leq
 c_5\mb^p_n(\bH)\delta_n^{p\rho(\alpha)},
\end{gather}
since $\omega(\alpha)\ge p \rho(\alpha)$ in both considered cases in view of the second equality in (\ref{eq-new-new:nu}) and of (\ref{eq:nu}).
Applying the third assertion of Theorem 3 (Part I), we obtain from
(\ref{eq12:preliminaries-proof-theorems}), (\ref{eq13:preliminaries-proof-theorems}), (\ref{eq96:preliminaries-proof-theorems}) and (\ref{eq100:preliminaries-proof-theorems})
\vskip0.1cm
\centerline{$
\displaystyle{\sup_{f\in\bN_{\vec{r},d}\big(\vec{\beta},\vec{L}\big)\cap\bF_{g}(R)}}\cR^{(p)}_n[\widehat{f}_{\vec{\mathbf{h}}(\cdot)}, f]\leq C
\bigg[(c_2+c_3+c_5+c_6)\mb^p_n(\bH)\delta_n^{p\rho(\alpha)}\bigg]^{\frac{1}{p}}
 \leq c_7\mb_n(\bH)\delta_n^{\rho(\alpha)},
$}
\vskip0.1cm
\noindent and the assertion of Theorem \ref{th:adaptive-L_p-bounded-deconv} follows in both considered cases.

\vskip0.1cm

$\mathbf{2^0.\;}$ Consider the case $\kappa_\alpha(p^*)<0, \tau(\infty)>0$.

Choose $\overline{\Bv}=\mathbf{v}$ and remark that the statements of Propositions \ref{prop1:auxiliary results-new} and \ref{prop2:auxiliary results-new} hold hold for any
$v\in[\underline{\Bv},\overline{\Bv}]$.
Indeed, $\tau(\infty)>0$ implies $\mathbf{v}<1$ and, therefore, $\overline{\mathbf{v}}=\overline{\mathbf{v}}\wedge\mathbf{v_3}=\mathbf{v}$.
We deduce from  (\ref{eq12:preliminaries-proof-theorems}), (\ref{eq13:preliminaries-proof-theorems}),
 (\ref{eq95:preliminaries-proof-theorems}) and (\ref{eq100:preliminaries-proof-theorems}), applying the first assertion of Theorem 3 (Part I)  that
 \vskip-0.7cm
\begin{gather}
\label{eq97:preliminaries-proof-theorems}
\sup_{f\in\bN_{\vec{r},d}\big(\vec{\beta},\vec{L}\big)\cap\bF_{g}(R)}\cR^{(p)}_n[\widehat{f}_{\vec{\mathbf{h}}(\cdot)}, f]\leq C
\bigg[c_8
 \delta_n^{\frac{\omega(\alpha)\tau(p)\beta(0)}{z(\alpha)}}+(c_3+c_6)\mb^p_n(\bH)\delta_n^{p\rho(\alpha)}\bigg]^{\frac{1}{p}}
 \leq c_9\mb_n(\bH)\delta_n^{\rho(\alpha)}.
\end{gather}
\vskip-0.2cm
\noindent Here we have also used (\ref{eq:nu}). This completes the proof of
Theorem \ref{th:adaptive-L_p-bounded-deconv}.

\subsection{Proof of Theorem \ref{th:adaptive-L_p-general-deconv}}

 In the following we assume $p^*<\infty$, since $p^*=\infty$ implies by definition of the anisotropic Nikol'skii class that $\bN_{\vec{r},d}\big(\vec{\beta},\vec{L}\big)\subset \bB_{\infty,d}(L_\infty)$. Hence, the results in that case follow from Theorem \ref{th:adaptive-L_p-bounded-deconv} since $\varrho(\alpha)=\rho(\alpha)$ when $p^*=\infty$.

 Moreover, we remark that the imposed condition $p>[\min_{j=1,\ldots} \mu_j]^{-1}$ implies $Y\geq[X+1]\mathbf{y}^{-1}-1/\mathbf{u}$ in view of (\ref{eq:rem:about-main-condition}) proved in Remark \ref{rem:after-prop1}. This, first, makes the second assertion of Proposition \ref{prop1:auxiliary results-new} applicable.

 Next, it allows (recall that $p^*<\infty$ and $\alpha=1$) to rewrite $\cI_\infty(1)$ appeared  in Proposition   \ref{prop2:auxiliary results-new} as

 \vskip0.15cm
 \centerline{$
 \cI_\infty(1)=[\Bv,\mathbf{v_3}]\mathrm{1}_{\{\kappa_1(p^*)\geq0\}}+[\Bv,\overline{\mathbf{v}}]\mathrm{1}_{\{\kappa_1(p^*)< 0\}}.
 $}

\vskip0.15cm

\vskip0.1cm

$\mathbf{1^0.\;}$ Consider the case $\kappa_\alpha(p^*)<0, \tau(p^*)>0$.

Taking
into account that $\vec{L}\in[L_0,L_\infty]$ we remark that in view of \cite{Nikolski} [Theorem 6.9.1, Section~6.9]
 $\bN_{\vec{r},d}\big(\vec{\beta},\vec{L}\big)\subset \bB_{p^*,d}(c_9L_\infty)$, where $c_9$ is independent of $\vec{L}$. Thus, Theorem 3 (Part I) is applicable with $\mathbf{u=\infty}$, $\mathbf{q}=p^*$ and $D=c_9L_\infty\vee Q$. Choose $\overline{\Bv}=\mathbf{v}$ and remark that the statements of Propositions \ref{prop1:auxiliary results-new} and \ref{prop2:auxiliary results-new} hold since $\overline{\mathbf{v}}=\mathbf{v}$. The assertion of the theorem is obtained from (\ref{eq12:preliminaries-proof-theorems}), (\ref{eq13:preliminaries-proof-theorems}),
 (\ref{eq95:preliminaries-proof-theorems}), (\ref{eq100:preliminaries-proof-theorems}), (\ref{eq-new:nu}) and the first assertion of Theorem 3 (Part I) by the same computations
 that led to (\ref{eq97:preliminaries-proof-theorems}).

\vskip0.1cm

$\mathbf{2^0.\;}$ Consider the case $\kappa_1(p^*)<0, \tau(p^*)\leq0$. Recall that   $p^*>p$ in this case because it is necessary for the existence of an uniformly consistent estimator.
Since the definition of the anisotropic Nikol'skii class implies that $\bN_{\vec{r},d}\big(\vec{\beta},\vec{L}\big)\subset \bB_{p^*,d}(L_\infty)$, we assert that
the second assertion of Theorem 3 (Part I) is applicable with  $\mathbf{u}=\infty$,  $\mathbf{q}=p^*$ and $D=L_\infty\vee Q$. Choose $\overline{\Bv}=\mathbf{v_2}$ and note that
$\overline{\mathbf{v}}=\mathbf{v_2}$ in the considered case. Thus, we deduce from (\ref{eq12:preliminaries-proof-theorems}), (\ref{eq13:preliminaries-proof-theorems}),
 (\ref{eq100:preliminaries-proof-theorems}) and (\ref{eq-new:nu})
$$
\sup_{f\in\bN_{\vec{r},d}\big(\vec{\beta},\vec{L}\big)\cap\bF_{g,\infty}(R,Q)}\cR^{(p)}_n[\widehat{f}_{\vec{\mathbf{h}}(\cdot)}, f]\leq C\bigg[c'_2\delta_n^{\omega(1)-\frac{\omega(1)\kappa_1(p,\mathbf{\infty})}{\kappa_1(p^*,\mathbf{\infty})}}+
(c_3+c_6)\mb^p_n(\bH)\delta_n^{p\varrho(\alpha)}
+ \delta_n^{\frac{\omega(1)(p-p^*)}{\kappa_1(p^*,\mathbf{\infty})}}\bigg]^{\frac{1}{p}},
$$
and the assertion of the theorem follows in this case.

\vskip0.1cm

$\mathbf{3^0.\;}$ It remains to study the case $\kappa_1(p^*)\geq0$.
Let $s$ be an arbitrary number satisfying (\ref{eq3:new-formulas}) of Lemma \ref{lem:auxilary-results-for-prpositions-1-and2}. Since $\tau(s)>0$ and $s>p^*$
we can assert in view of \cite{Nikolski} [Theorem 6.9.1, Section~6.9]
 $\bN_{\vec{r},d}\big(\vec{\beta},\vec{L}\big)\subset \bB_{s,d}(c_{9}L_\infty)$, where $c_{9}$ is independent of $\vec{L}$. Thus, Theorem 3 (Part) is applicable with $\mathbf{u=\infty}$, $\mathbf{q}=s$ and $D=c_{10}L_\infty\vee Q$.
 Choosing $\overline{\Bv}=\mathbf{v_3}$, we deduce from (\ref{eq12:preliminaries-proof-theorems}), (\ref{eq13:preliminaries-proof-theorems}),
 (\ref{eq100:preliminaries-proof-theorems}) and from the second assertion of Theorem 3 (Part 2)
\begin{gather}
\label{eq98:preliminaries-proof-theorems}
\sup_{f\in\bN_{\vec{r},d}\big(\vec{\beta},\vec{L}\big)\cap\bF_{g,\infty}(R,Q)}\cR^{(p)}_n[\widehat{f}_{\vec{\mathbf{h}}(\cdot)}, f]\leq C\bigg[(c_3+c_6)\mb^p_n(\bH)\delta_n^{p\varrho(\alpha)}
+\delta_n^{\frac{Y(s-p)}{[1+X]/\omega(0)-Y/\beta(0)}}\bigg]^{\frac{1}{p}}.
\end{gather}
Since either  $p^*/\omega(0)=1/\beta(0)$, $\kappa_{1}(p^*)> 0$ or $p^*/\omega(0)>1/\beta(0)$, $\kappa_1(p^*)\geq 0$ and $s>p^*\geq p $ we get
$$
\frac{Y(s-p)}{[1+X]/\omega(0)-Y/\beta(0)}=\frac{[\omega(0)-\omega(1)](s-p)}{\kappa_1(p^*,\infty)+p^*-\omega(0)/\beta(0)}>0
$$
Simple algebra shows that
$$
\frac{Y(s-p)}{[1+X]/\omega(0)-Y/\beta(0)}\geq \frac{p}{2+1/\beta(1)}\quad\Leftrightarrow\quad
sY\omega(1)\geq p(X+1)[2+1/\beta(1)]^{-1}
$$
Using again $\kappa_1(p^*)\geq 0$ and $p^*\geq p$ we obtain
\begin{eqnarray*}
sY\omega(1)\big[2+1/\beta(1)\big]&=&s\kappa_1(p^*)Y+sp^*Y\geq spY\geq p(X+1)
\end{eqnarray*}
since  $s$ satisfies  (\ref{eq3:new-formulas}) of Lemma \ref{lem:auxilary-results-for-prpositions-1-and2}. Thus, we have for all $n$ large enough
\vskip0.2cm
\centerline{$
\delta_n^{\frac{Y(s-p)}{[1+X]/\omega(0)-Y/\beta(0)}}\leq\delta_n^{\frac{p}{2+1/\beta(1)}}\leq\mb^p_n(\bH)\delta_n^{p\varrho(\alpha)}
$}
\vskip0.2cm
\noindent and the assertion of the theorem in the case $\kappa_1(p^*)\geq 0$ follows from  (\ref{eq98:preliminaries-proof-theorems}) and the first equality in (\ref{eq-new-new:nu}).
Theorem \ref{th:adaptive-L_p-general-deconv} is proved.

\section{Proofs of Propositions \ref{prop1:auxiliary results-new}, \ref{prop2:auxiliary results-new} and \ref{prop:measure-of-bias-deconv}}
\label{sec:Proofs-of-Propositions} The proof of Lemma \ref{lem:relations-between-parameters} is postponed to  Appendix.

\begin{lemma}
\label{lem:relations-between-parameters} For any $\vec{\beta}$, $\vec{r}$, $\vec{\mu}$, $p\geq 1$ and $\alpha\in[0,1]$ the following is true.

\vskip0.2cm

\centerline{$
1/\gamma(\alpha)-1/\beta(\alpha)=\big[\tau(\infty)\beta(0)\big]^{-1}
\big[1/\omega(\alpha)-1/\upsilon(\alpha)\big].
$}

\end{lemma}

\subsection{Proof of Proposition \ref{prop1:auxiliary results-new}}

We start the proof with several remarks which will be useful in the sequel.
First, obviously  there exists $0<\mathbf{T}:=T\big(\vec{\beta},\vec{r},\vec{\mu},p\big)<\infty$ independent of $\vec{L}$ such that
\begin{equation}
\label{eq01:proof-prop1:auxiliary results-new}
\lim_{n\to\infty}\;(\ln{n})^{-1}\sup_{\alpha\in\{0,1\}}\;\sup_{\mathbf{s}\in[1,\infty]}
\sup_{v\in[\underline{\mathbf{v}},1\vee\mathbf{v}]}\sum_{j=1}^d\Big\{\big|\ln{\big(\blh_j(v,\mathbf{1})\big)}\big|+\big|\ln{\big(\mh_j(v,\mathbf{s})\big)}\big|\Big\}= \mathbf{T}.
\end{equation}
Next, for any $\mathbf{s}\in [1,\infty]$ and any $v>0$
\begin{equation}
\label{eq02:proof-prop1:auxiliary results-new}
\frac{\ln{n}}{n}\prod_{j=1}^d\big(\widetilde{\bleta}_j(v,\mathbf{s})\big)^{-1-2\blg_j(\alpha)}=\ma^2\BL^{-\frac{1}{\beta(\alpha)}}
\big(\ma^{-2}\delta_n\big)^{\frac{\omega(\alpha)}{\omega(\alpha)+\mathbf{s}}}
v^{\frac{2\mathbf{s}-\omega(\alpha)/\beta(\alpha)}{\mathbf{s}+\omega(\alpha)}}=2\BL^{-\frac{1}{\beta(\alpha)}}\ma^2v^2\mz^{-1}(v).
\end{equation}

1) Let us proceed to the  proof of the first assertion. First we remark that for all $n\geq 3$
\begin{equation}
\label{eq1:proof-prop1:auxiliary results-new}
\vec{\blh}(v,\mathbf{1})\in (0,1]^d,\quad \forall v\in[\underline{\mathbf{v}},1].
\end{equation}
Indeed for any $v>0$ we have
 since  $\BL\leq L_0$,
\begin{equation}
\label{eq2:proof-prop1:auxiliary results-new}
\widetilde{\bleta}^{\beta_jr_j}_j(v,\mathbf{1})\leq \big(\ma^{-2}\delta_n\big)^{\frac{\omega(\alpha)}{\mathbf{1}+\omega(\alpha)}}
v^{r_j-\frac{\omega(\alpha)(2+1/\beta(\alpha))}{\mathbf{1}+\omega(\alpha)}}, \quad j\in\bar{\cJ}_\infty.
\end{equation}
Therefore, for any $v\in[\underline{\mathbf{v}},1]$ one has in view of the definition of $\underline{\mathbf{v}}$
$$
\widetilde{\bleta}^{\beta_jr_j}_j(v,\mathbf{1})\leq \big(\ma^{-2}\delta_n\big)^{\frac{\omega(\alpha)}{\mathbf{1}+\omega(\alpha)}}
v^{1-\frac{\omega(\alpha)(2+1/\beta(\alpha))}{\mathbf{1}+\omega(\alpha)}}\leq \big(\ma^{-2}\delta_n\big)^{\frac{\omega(\alpha)}{\mathbf{1}+\omega(\alpha)}}
\underline{\mathbf{v}}^{1-\frac{\omega(\alpha)(2+1/\beta(\alpha))}{\mathbf{1}+\omega(\alpha)}}=1, \quad j\in\bar{\cJ}_\infty.
$$
Note that for any  $j\in\cJ_\infty$
$$
\widetilde{\bleta}_j(v,\mathbf{1})=\big(\BL L_j^{-1}v\big)^{\frac{1}{\beta_j}}\leq v^{\frac{1}{\beta_j}}\leq 1,\quad \quad \forall v\leq 1.
$$
and the proof of (\ref{eq1:proof-prop1:auxiliary results-new}) is completed since $\blh_j(\cdot,\mathbf{1})\leq \widetilde{\bleta}_j(\cdot,\mathbf{1})$ by construction.

Set $T_0=\big[\mathbf{T}+2\big]\; e^{d+2\sum_{j=1}^d\mu_j(\alpha)}\BL^{-\frac{1}{\beta(\alpha)}}$ and
remark that in view of (\ref{eq01:proof-prop1:auxiliary results-new}), (\ref{eq02:proof-prop1:auxiliary results-new}) and (\ref{eq1:proof-prop1:auxiliary results-new}) for all $n$ large enough and any $v\in[\underline{\mathbf{v}},1]$
\begin{eqnarray}
\label{eq3:proof-prop1:auxiliary results-new}
 G_n\big(\vec{\blh}(v,\mathbf{1})\big)&\leq&\frac{(\mathbf{T}+2)\ln{n}}{n\prod_{j=1}^d\big(\blh_j(v,\mathbf{1})\big)^{1+\blg_j(\alpha)}}
\leq \frac{T_0\BL^{\frac{1}{\beta(\alpha)}}\ln{n}}{n\prod_{j=1}^d\big(\widetilde{\bleta}_j(v,\mathbf{1})\big)^{1+\blg_j(\alpha)}}
\nonumber\\
&\leq& \frac{T_0\BL^{\frac{1}{\beta(\alpha)}}\ln{n}}{n\prod_{j=1}^d\big(\widetilde{\bleta}_j(v,\mathbf{1})\big)^{1+2\blg_j(\alpha)}}
=T_0\ma^{\frac{2}{1+\omega(\alpha)}}\delta_n^{\frac{\omega(\alpha)}{1+\omega(\alpha)}}
v^{\frac{2-\omega(\alpha)/\beta(\alpha)}{\mathbf{1}+\omega(\alpha)}}.
\end{eqnarray}
Here we have taken into account that $\blh_j(v,\mathbf{s})\geq e^{-1}\bleta_j(v,\mathbf{s})$. Since
$$
T_0\ma^{\frac{2}{1+\omega(\alpha)}}\delta_n^{\frac{\omega(\alpha)}{\omega(\alpha)+1}}
v^{\frac{2-\omega(\alpha)/\beta(\alpha)}{1+\omega(\alpha)}}\leq T_0\ma^2 v\;\;
\Leftrightarrow\;\; v\geq \underline{\mathbf{v}},
$$
denoting $\ma=\sqrt{a/T_0}$ we assert that
$$
G_n\big(\vec{\blh}(v,\mathbf{1})\big)\leq a v, \quad \forall v\in [\underline{\mathbf{v}},1].
$$
The first assertion is established.

2)  Before proving the second assertion, let us make several remarks.

\vskip0.1cm

$\mathbf{1^0.}\;$ For any $\mathbf{u}\in[1,\infty]$ the following is true.
\begin{equation}
\label{eq50:proof-prop1:auxiliary results-new}
\widehat{\bleta}_j(\mathbf{v},\mathbf{u})=\big(\BL L_j^{-1}\mathbf{v}\big)^{\frac{1}{\beta_j}},\; j\in\cJ_\infty,\quad\;
\widehat{\bleta}_j(\mathbf{v},\mathbf{u})=\big(\BL L_j^{-1}\big)^{\frac{1}{\gamma_j}} \big(\ma^{-2}\delta_n\big)^{\frac{\omega(\alpha)\tau(p_\pm)\beta(0)}{\gamma_j[z(\alpha)+\omega(\alpha)/\mathbf{u}]}}
,\; j\in\bar{\cJ}_\infty.
\end{equation}
The first equality follows directly from the definition of $\widehat{\bleta}_j(\mathbf{v},\mathbf{u})$ since, remind $\gamma_j=\beta_j, q_j=\infty$ if $j\in\cJ_\infty$. Thus, let us prove the second equality. We have
$$
\widehat{\bleta}^{\gamma_jq_j}_j(\mathbf{v},\mathbf{u})=\big(\BL L_j^{-1}\big)^{p_\pm}
\big(\ma^{-2}\delta_n\big)^{\frac{\mathbf{u}\upsilon(\alpha)}{\mathbf{u}+\upsilon(\alpha)}}
 \mathbf{v}^{p_\pm-\frac{\mathbf{u}\upsilon(\alpha)(2+1/\gamma(\alpha))}{\mathbf{u}+\upsilon(\alpha)}},\quad \forall j\in\bar{\cJ}_\infty.
$$
Here we have used that $q_j=p_\pm$ for any $j\in\bar{\cJ}_\infty$. Using the definition of $\mathbf{v}$ we get
$$
\widehat{\bleta}^{\gamma_jq_j}_j(\mathbf{v},\mathbf{u})=\big(\BL L_j^{-1}\big)^{p_\pm} \big(\ma^{-2}\delta_n\big)^{\frac{\mathbf{u}\upsilon(\alpha)}{\mathbf{u}+\upsilon(\alpha)}+
\frac{\omega(\alpha)\tau(\infty)\beta(0)}{z(\alpha)+\omega(\alpha)/\mathbf{u}}
\left[p_\pm-\frac{\mathbf{u}\upsilon(\alpha)(2+1/\gamma(\alpha))}{\mathbf{u}+\upsilon(\alpha)}\right]}
,\quad \forall j\in\bar{\cJ}_\infty.
$$
Using the definition of $z(\alpha)$ we obtain

\begin{eqnarray*}
A&:=&\frac{\mathbf{u}\upsilon(\alpha)}{\mathbf{u}+\upsilon(\alpha)}+
\frac{\omega(\alpha)\tau(\infty)\beta(0)}{z(\alpha)+\omega(\alpha)/\mathbf{u}}
\left[p_\pm-\frac{\mathbf{u}\upsilon(\alpha)(2+1/\gamma(\alpha))}{\mathbf{u}+\upsilon(\alpha)}\right]
=\frac{\omega(\alpha)\tau(\infty)\beta(0)p_\pm}{z(\alpha)+\omega(\alpha)/\mathbf{u}}
\\
&&+\frac{\mathbf{u}\upsilon(\alpha)}{(\mathbf{u}+\upsilon(\alpha))(z(\alpha)+\omega(\alpha)/\mathbf{u})}
\big[1+\omega(\alpha)/\mathbf{u}-\omega(\alpha)\tau(\infty)\beta(0)\big\{1/\gamma(\alpha)
-1/\beta(\alpha)\big\}\big].
\end{eqnarray*}
We obtain applying  Lemma \ref{lem:relations-between-parameters}
$$
A=\frac{\omega(\alpha)\tau(\infty)\beta(0)p_\pm}{z(\alpha)+\omega(\alpha)/\mathbf{u}}+
\frac{\mathbf{u}\upsilon(\alpha)\omega(\alpha)[1/s+1/\upsilon(\alpha)]}{(\mathbf{u}+\upsilon(\alpha))(z(\alpha)+\omega(\alpha)/\mathbf{u})}=
\frac{\omega(\alpha)\tau(p_\pm)\beta(0)p_\pm}{z(\alpha)+\omega(\alpha)/\mathbf{u}}.
$$
The second formula in (\ref{eq50:proof-prop1:auxiliary results-new}) is established.

$\mathbf{2^0.}\;$ Next, let us prove that
\begin{equation}
\label{eq5:proof-prop1:auxiliary results-new}
\vec{\mh}(\mathbf{v},\mathbf{u})\in (0,1]^d, \quad \forall \mathbf{u}\in [1,\infty].
\end{equation}
If $\cJ_\infty\neq \emptyset$, which is equivalent to $p^*=\infty$, the definition of $\mathbf{v}$ implies that $\mathbf{v}\leq 1$ for all $n$ large enough,
since $\tau(p^*)=\tau(\infty)> 0$ and in view of (\ref{eq4:proof-prop1:auxiliary results-new}). We deduce from the first equality in (\ref{eq50:proof-prop1:auxiliary results-new})
\begin{equation*}
\mh_j(\mathbf{v},\mathbf{u})\leq\widehat{\bleta}_j(\mathbf{v},\mathbf{u})=\big(\BL L_j^{-1}\mathbf{v}\big)^{\frac{1}{\beta_j}}\leq \mathbf{v}^{\frac{1}{\beta_j}}\leq 1,\quad \forall j\in\cJ_\infty.
\end{equation*}
and (\ref{eq5:proof-prop1:auxiliary results-new}) is proved for any $j\in\cJ_\infty$.

It remains to note that $\tau(p_\pm)\geq \tau(p^*)$ since $p^*\geq p_\pm$ and therefore, if $\tau(p^*)\geq 0$ we have
$$
\mh_j(\mathbf{v},\mathbf{u})\leq\widehat{\bleta}_j(\mathbf{v},\mathbf{u})\leq 1
,\quad \forall j\in\bar{\cJ}_\infty,
$$
for all $n$ large enough in view of (\ref{eq4:proof-prop1:auxiliary results-new}) of Lemma \ref{lem:auxilary-results-for-prpositions-1-and2}, the second equality in (\ref{eq50:proof-prop1:auxiliary results-new})  and since $\BL L_j^{-1}\leq 1$. This completes the proof of (\ref{eq5:proof-prop1:auxiliary results-new}).

\vskip0.1cm

$\mathbf{3^0.}\;$  For any $\mathbf{u}\in[1,\infty]$ one has
\begin{eqnarray}
\label{eq500:proof-prop1:auxiliary results-new}
\ma^{-2}\delta_n\prod_{j=1}^d\widehat{\bleta}^{-1-2\blg_j(\alpha)}_j(\mathbf{v},\mathbf{u})
&\leq&T^{-1}(\alpha)\big(\ma^{-2}\delta_n\big)^{1-\frac{\omega(\alpha)\tau(\infty)\beta(0)/\beta(\alpha)+1}{z(\alpha)+\omega(\alpha)/\mathbf{u}}};
\\
\label{eq501:proof-prop1:auxiliary results-new}
\ma^{-2}\delta_n\prod_{j=1}^d\widehat{\bleta}^{-1}_j(\mathbf{v},
\mathbf{u})&\leq&T^{-1}(0)\big(\ma^{-2}\delta_n\big)^{1-\frac{\omega(\alpha)}{z(\alpha)+\omega(\alpha)/\mathbf{u}}}
\end{eqnarray}
where we have denoted  $T\big(\alpha\big)=\inf_{\vec{L}\in[L_0,L_\infty]^d}\prod_{j\in\cJ_\infty}\big(\BL L_j^{-1}\big)^{\frac{1+2\blg_j(\alpha)}{\beta_j}}\prod_{j\in\bar{\cJ}_\infty}\big(\BL L_j^{-1}\big)^{\frac{1+2\blg_j(\alpha)}{\gamma_j}}$.

Indeed, we have in view of (\ref{eq50:proof-prop1:auxiliary results-new}) and the definition of $\mathbf{v}$
\begin{eqnarray*}
\label{eq8:proof-prop1:auxiliary results-new}
\prod_{j=1}^d\widehat{\bleta}^{1+2\blg_j(\alpha)}_j(\mathbf{v},\mathbf{u})&\geq&T^{-1}(\alpha)
\big(\ma^{-2}\delta_n\big)^{\frac{\omega(\alpha)\tau(p_\pm)\beta(0)}{\gamma_\pm(\alpha)[z(\alpha)+\omega(\alpha)/\mathbf{u}]}
+\frac{\omega(\alpha)\tau(\infty)\beta(0)}{\beta_\infty(\alpha)[z(\alpha)+\omega(\alpha)/\mathbf{u}]}}
\\
\prod_{j=1}^d\widehat{\bleta}_j(\mathbf{v},\mathbf{u})&\geq&T^{-1}(0)
\big(\ma^{-2}\delta_n\big)^{\frac{\omega(\alpha)\tau(p_\pm)\beta(0)}{\gamma_\pm(0)[z(\alpha)+\omega(\alpha)/\mathbf{u}]}
+\frac{\omega(\alpha)\tau(\infty)\beta(0)}{\beta_\infty(0)[z(\alpha)+\omega(\alpha)/\mathbf{u}]}},
\end{eqnarray*}
where we have put
$
\frac{1}{\beta_\infty(\alpha)}=\sum_{j\in\cJ_\infty}\frac{1+2\blg(\alpha)}{\beta_j},\;
\frac{1}{\gamma_\pm(\alpha)}=
\sum_{j\in\bar{\cJ}_\infty}\frac{1+2\blg(\alpha)}{\gamma_j}.
$
Note that for any $\alpha\in[0,1]$
\begin{eqnarray*}
\frac{\tau(p_\pm)}{\gamma_\pm(\alpha)}+\frac{\tau(\infty)}{\beta_\infty(\alpha)}=\sum_{j\in\bar{\cJ}_\infty}
\frac{(1+2\blg(\alpha))\tau(r_j)}{\beta_j}+\frac{\tau(\infty)}{\beta_\infty(\alpha)}
=\frac{\tau(\infty)}{\beta(\alpha)}+\frac{1}{\omega(\alpha)\beta(0)}.
\end{eqnarray*}
and (\ref{eq500:proof-prop1:auxiliary results-new}) and (\ref{eq501:proof-prop1:auxiliary results-new}) are established.

\vskip0.1cm

$\mathbf{4^0.}\;$ Simple algebra shows that for any $\mathbf{u}\in[1,\infty]$
$$
\big(\ma^{-2}\delta_n\big)^{1-\frac{\omega(\alpha)\tau(\infty)\beta(0)/\beta(\alpha)+1}{z(\alpha)+\omega(\alpha)/\mathbf{u}}}=2\mathbf{v}^2\mz^{-1}(\mathbf{v}),
$$
and we deduce from (\ref{eq500:proof-prop1:auxiliary results-new}) for any $\mathbf{u}\in[1,\infty]$ (recall that $\mz\equiv 2$ if $\mathbf{u}=\infty$)
\begin{eqnarray}
\label{eq502:proof-prop1:auxiliary results-new}
&&\delta_n\prod_{j=1}^d\widehat{\bleta}^{-1-2\blg_j(\alpha)}_j(\mathbf{v},\mathbf{u})
\leq 2T^{-1}(\alpha)\ma^2\mathbf{v}^2\mz^{-1}(\mathbf{v}).
\end{eqnarray}
Let us also prove that for any $\mathbf{u}\in[1,\infty]$ and all $n$ large enough
\begin{eqnarray}
\label{eq503:proof-prop1:auxiliary results-new}
&&\mathbf{v}>\Bv:=\big(\ma^{-2}\delta_n\big)^{\frac{1}{2+1/\beta(\alpha)}}\quad \Rightarrow\quad \mz(\mathbf{v})\geq 2.
\end{eqnarray}
The latter inclusion follows from (\ref{eq02:new-paper}).
Indeed, if $\tau(\infty)\leq 0$ then $\mathbf{v}\geq 1 \ge \Bv$. If $\tau(\infty)>0$
$$
\frac{\omega(\alpha)\tau(\infty)\beta(0)}{z(\alpha)+\omega(\alpha)/\mathbf{u}}-\frac{1}{2+1/\beta(\alpha)}=
-\frac{1+\omega(\alpha)/\mathbf{u}}{[z(\alpha)+\omega(\alpha)/\mathbf{u}][2+1/\beta(\alpha)]}<0
$$
in view of (\ref{eq4:proof-prop1:auxiliary results-new}), so $\mathbf{v}>\Bv$.
Note at last that for any $\mathbf{u}\in [1,\infty]$
\begin{eqnarray}
\label{eq504:proof-prop1:auxiliary results-new}
\mathbf{v}\mz^{-1}(\mathbf{v})=2\big(\ma^{-2}\delta_n\big)^{\frac{\omega(\alpha)\tau(\mathbf{u})\beta(0)}{z(\alpha)+\omega(\alpha)/\mathbf{u}}}.
\end{eqnarray}

$\mathbf{5^0.}\;$ Let us proceed to the proof of the second assertion. Let us choose $\ma<aT(\alpha)/(4T_0)<1$. We have in view of (\ref{eq01:proof-prop1:auxiliary results-new}), (\ref{eq500:proof-prop1:auxiliary results-new}) and (\ref{eq502:proof-prop1:auxiliary results-new}) similarly to (\ref{eq3:proof-prop1:auxiliary results-new})
\begin{eqnarray}
\label{eq507:proof-prop1:auxiliary results-new}
&& F^2_n\big(\vec{\mh}(\mathbf{v},\mathbf{u})\big)\leq\frac{T_0\delta_n}{\prod_{j=1}^d
 \big(\widehat{\bleta}_j(\mathbf{v},\mathbf{u})\big)^{1+2\blg_j(\alpha)}}
\leq 2T_0T^{-1}(\alpha)\ma^2\mathbf{v}^2\mz^{-1}(\mathbf{v})\leq a^2\mathbf{v}^2\mz^{-1}(\mathbf{v}).
\end{eqnarray}

Thus to prove the assertion all we need to show is that $\vec{\mh}(\mathbf{v},\mathbf{u})\in\mH(\mathbf{v})$, i.e. $G_n\big(\vec{\mh}(\mathbf{v},\mathbf{u})\big)\le a\mathbf{v}$. Let us distinguish three cases.

\vskip0.1cm

$\mathbf{5^0a.}\;$ Let $\tau(\infty)\geq 0.$ We remark that the definition of $\mathbf{v}$ in this case yields $\mathbf{v}\leq 1$ for all $n$ large enough
and we obtain from (\ref{eq502:proof-prop1:auxiliary results-new}) and (\ref{eq503:proof-prop1:auxiliary results-new}) that
\begin{eqnarray}
\label{eq505:proof-prop1:auxiliary results-new}
&&\delta_n\prod_{j=1}^d\widehat{\bleta}^{-1-2\blg_j(\alpha)}_j(\mathbf{v},\mathbf{u})
\leq T^{-1}(\alpha)\ma^2\mathbf{v}.
\end{eqnarray}
Then we have in view of (\ref{eq01:proof-prop1:auxiliary results-new}), (\ref{eq5:proof-prop1:auxiliary results-new}), (\ref{eq500:proof-prop1:auxiliary results-new}) and
(\ref{eq505:proof-prop1:auxiliary results-new}) similarly to (\ref{eq3:proof-prop1:auxiliary results-new})
\begin{eqnarray}
\label{eq506:proof-prop1:auxiliary results-new}
 &&G_n\big(\vec{\mh}(\mathbf{v},\mathbf{u})\big)\leq
\frac{T_0\delta_n}{\prod_{j=1}^d
 \big(\widehat{\bleta}_j(\mathbf{v},\mathbf{u})\big)^{1+2\blg_j(\alpha)}}
\leq T_0T^{-1}(\alpha)\ma^2\mathbf{v}\le a\mathbf{v}.
\end{eqnarray}

\vskip0.1cm

$\mathbf{5^0b.}\;$ Let $\tau(\infty)< 0, \tau(p^*)> 0$ and $\alpha\neq 1$. Then by assumption $\mathbf{u}\le p$, and thus $\tau (\mathbf{u})\ge 0$. We get from (\ref{eq502:proof-prop1:auxiliary results-new}) and (\ref{eq504:proof-prop1:auxiliary results-new})
\begin{eqnarray}
\label{eq508:proof-prop1:auxiliary results-new}
&&\delta_n\prod_{j=1}^d\widehat{\bleta}^{-1-2\blg_j(\alpha)}_j(\mathbf{v},\mathbf{u}^*)
\leq 4T^{-1}(\alpha)\ma^2\mathbf{v},
\end{eqnarray}
so $G_n\big(\vec{\mh}(\mathbf{v},\mathbf{u})\big)\le a\mathbf{v}$ follows from (\ref{eq506:proof-prop1:auxiliary results-new}) and (\ref{eq507:proof-prop1:auxiliary results-new}).

\vskip0.1cm

$\mathbf{5^0c.}\;$ Let $\tau(\infty)< 0, \tau(p^*)> 0$, $\alpha=1$. We have as previously
\begin{eqnarray}
\label{eq6:proof-prop1:auxiliary results-new}
G^2_n\big(\vec{h}(\mathbf{u})\big)&\leq& \frac{(\mathbf{T}+2)\ln{n}}{n\prod_{j=1}^d\big(\mh_j(\mathbf{v},\mathbf{u})\big)
^{1+2\blg_j(\alpha)}}\frac{(\mathbf{T}+2)\ln{n}}{n\prod_{j=1}^d \mh_j(\mathbf{v},\mathbf{u})}
\nonumber\\
&\leq&
2T^2_0T^{-1}(1)T_1\ma^4\mathbf{v}^2\mz^{-1}(\mathbf{v})\bigg[\frac{T(0)\ma^{-2}\delta_n}{\prod_{j=1}^d\widehat{\bleta}_j(\mathbf{v},\mathbf{u})}\bigg].
\end{eqnarray}
Here we have used (\ref{eq502:proof-prop1:auxiliary results-new}) and put $T_1=T^{-1}(0)\BL^{-1/\beta(0)}$.
Our goal now is to show that for any $\mathbf{u}\in [1,\infty]$ and all $n$ large enough
\begin{equation}
\label{eq7:proof-prop1:auxiliary results-new}
T(0)\ma^{-2}\delta_n\mz^{-1}(\mathbf{v})\prod_{j=1}^d\widehat{\bleta}^{-1}_j(\mathbf{v},\mathbf{u})\leq 1.
\end{equation}
In view of  (\ref{eq501:proof-prop1:auxiliary results-new}) and of the definition of $\mz(\cdot)$
 in order to establish  (\ref{eq7:proof-prop1:auxiliary results-new}) it suffices to show that
$ z(1)/\omega(1)-1+2/\mathbf{u}\geq 0$.

Since we assumed $\tau(\infty)<0$ and $\tau(p^*)>0$, then necessarily $\mathbf{u}^*>p^*$ since $\tau(\mathbf{u}^*)=0$ and $\tau(\cdot)$ is strictly decreasing. Hence, the required results follows from (\ref{eq900:proof-prop1:auxiliary results-new}). Thus, (\ref{eq7:proof-prop1:auxiliary results-new}) is proved. Then choosing $\ma$ such that $T_0(2T^{-1}(1)T_1)^{1/2}\ma^2\leq a$, we obtain
from (\ref{eq6:proof-prop1:auxiliary results-new}) and (\ref{eq7:proof-prop1:auxiliary results-new}) that for all all
$n$ large enough
$$
G_n\big(\vec{\mh}(\mathbf{v},\mathbf{u})\big)\leq T_0(2T^{-1}(\alpha)T_1)^{1/2}\ma^2\mathbf{v}\leq av.
$$
The second assertion is proved
\epr


\subsection{Proof of Proposition \ref{prop2:auxiliary results-new}}
We start the proof with several remarks which will be useful in the sequel.

$\mathbf{1^0.}\;$
Let us show that  for all $n$ large enough
\begin{equation}
\label{eq1:proof-prop2:auxiliary results-new}
\vec{\blh}(v,\mathbf{u})\in (0,1]^d, \quad   \forall  v\in\cI_{\mathbf{u}}(\alpha),\;\;\forall\mathbf{u}\geq\mathbf{u}^*\vee p^*.
\end{equation}
In view of the definition of $\widetilde{\bleta}_j(\cdot,\mathbf{u}), j=1,\ldots,d$,
\begin{equation}
\label{eq2:proof-prop2:auxiliary results-new}
\widetilde{\bleta}^{\beta_jr_j}_j(v,\mathbf{u})= \big(\BL L_j^{-1}\big)^{r_j}\big\{\ma^{-2}\delta_n\big\}^{\frac{\mathbf{u}\omega(\alpha)}{\mathbf{u}+\omega(\alpha)}}
v^{r_j-\frac{\mathbf{u}\omega(\alpha)(2+1/\beta(\alpha))}{\mathbf{u}+\omega(\alpha)}}, \quad j\in\bar{\cJ}_\infty.
\end{equation}
Therefore, for any $v\in[\mathbf{v}_0,1]$ one has, taking into account that $\BL\leq L_0$,
$$
\widetilde{\bleta}^{\beta_jr_j}_j(v,\mathbf{u})\leq \big\{\ma^{-2}\delta_n\big\}^{\frac{\mathbf{u}\omega(\alpha)}{\mathbf{u}+\omega(\alpha)}}
v^{1-\frac{\mathbf{u}\omega(\alpha)(2+1/\beta(\alpha))}{\mathbf{u}+\omega(\alpha)}}\leq \big\{\ma^{-2}\delta_n\big\}^{\frac{\mathbf{u}\omega(\alpha)}{\mathbf{u}+\omega(\alpha)}}
\mathbf{v}_0^{1-\frac{\mathbf{u}\omega(\alpha)(2+1/\beta(\alpha))}{\mathbf{u}+\omega(\alpha)}}=1, \quad j\in\bar{\cJ}_\infty.
$$
It remains to note that $\Bv>\mathbf{v}_0$ for all $n$ large enough and, therefore,
\begin{equation}
\label{eq3:proof-prop2:auxiliary results-new}
\widetilde{\bleta}_j(v,\mathbf{u}) \le 1,\;j\in\bar{\cJ}_\infty,\quad\forall v\in[\Bv,1]\cap\cI_\mathbf{u}(\alpha).
\end{equation}
We also have in view of the definition of $\widetilde{\bleta}_j(\cdot,\mathbf{u}), j=1,\ldots,d$,
\begin{equation*}
\widetilde{\bleta}_j(v,\mathbf{u})= \big(\BL L_j^{-1}v\big)^{\frac{1}{\beta_j}}\leq 1, \quad j\in\cJ_\infty,
\end{equation*}
for any $v\leq 1$. This together with (\ref{eq3:proof-prop2:auxiliary results-new}) proves (\ref{eq1:proof-prop2:auxiliary results-new}) in the cases when $\cI_{\mathbf{u}}(\alpha)=[\Bv,1]$.

Noting that $p^*<\infty$ is equivalent to $\cJ_\infty=\emptyset$, we deduce from (\ref{eq2:proof-prop2:auxiliary results-new})
for any $v\geq1$
$$
\widetilde{\bleta}^{\beta_jr_j}_j(v,\mathbf{u})\leq \big\{\ma^{-2}\delta_n\big\}^{\frac{\mathbf{u}\omega(\alpha)}{\mathbf{u}+\omega(\alpha)}}
v^{p^*-\frac{\mathbf{u}\omega(\alpha)(2+1/\beta(\alpha))}{\mathbf{u}+\omega(\alpha)}}\leq \big\{\ma^{-2}\delta_n\big\}^{\frac{\mathbf{u}\omega(\alpha)}{\mathbf{u}+\omega(\alpha)}}
v^{-\kappa_\alpha(p^*,\mathbf{u})}, \; j=1,\ldots,d.
$$
Thus, if $\kappa_\alpha(p^*,\mathbf{u})\geq 0$ then for any $v\geq 1$
$$
\widetilde{\bleta}^{\beta_jr_j}_j(v,\mathbf{u})\leq \big\{\ma^{-2}\delta_n\big\}^{\frac{\mathbf{u}\omega(\alpha)}{\mathbf{u}+\omega(\alpha)}}\to 0, n\to\infty, \;\; j=1,\ldots,d.
$$
This together with (\ref{eq3:proof-prop2:auxiliary results-new}) yields (\ref{eq1:proof-prop2:auxiliary results-new}) in the case $\kappa_\alpha(p^*,\mathbf{u})\geq 0, p^*<\infty$, whatever the value of $\alpha$.

Let  $\alpha=1, p^*<\infty, \kappa_\alpha(p^*,\mathbf{u})< 0, \tau(p^*)>0$.

\noindent Then $\overline{\mathbf{v}}=\mathbf{v}$  and we have for any $j=1,\ldots,d$ and $v\in [1,\mathbf{v}]$ in view of the definition of $\mathbf{v}$
$$
\widetilde{\bleta}^{\beta_jr_j}_j(v,\mathbf{u})\leq\widetilde{\bleta}^{\beta_jr_j}_j(\mathbf{v},\mathbf{u})= \big\{\ma^{-2}\delta_n\big\}^{\frac{\mathbf{u}\omega(1)}{\mathbf{u}+\omega(1)}
-\frac{\kappa_1(p^*,\mathbf{u})\omega(1)\tau(\infty)\beta(0)}{z(1)+\omega(1)/\mathbf{u}}}=
\big\{\ma^{-2}\delta_n\big\}^{\frac{p^*\tau(p^*)\omega(1)}{z(1)+\omega(1)/\mathbf{u}}}\to 0,\;n\to\infty,
$$
in view of (\ref{eq4:proof-prop1:auxiliary results-new}). Hence, (\ref{eq1:proof-prop2:auxiliary results-new}) holds in this case.

Let $\alpha=1, \kappa_\alpha(p^*,\mathbf{u})< 0, \tau(p^*)\leq 0$.

\noindent Then $\overline{\mathbf{v}}=\mathbf{v}_2$  and we have for any $v\in [1,\mathbf{v_2}]$ in view of the definition of $\mathbf{v_2}$
$$
\widetilde{\bleta}^{\beta_jr_j}_j(v,\mathbf{u})\leq \widetilde{\bleta}^{\beta_jr_j}_j(\mathbf{v_2},\mathbf{u})=\big\{\ma^{-2}\delta_n\big\}^{\frac{\mathbf{u}\omega(1)}{\mathbf{u}+\omega(1)}}
\mathbf{v}_2^{-\kappa_1(p^*,\mathbf{u})}=1, \quad j=1,\ldots,d.
$$
and, therefore (\ref{eq1:proof-prop2:auxiliary results-new}) holds in this case.

Let $\kappa_\alpha(p^*,\mathbf{u})< 0, \;\alpha\neq 1, \mathbf{u}<\infty$. First we note that $\tau(\infty)<0$ and  $\mathbf{u}\geq \mathbf{u^*}\vee p^*$ imply
$$
1-\mathbf{u}/\omega(0)+1/\beta(0)=1-\mathbf{u}+\mathbf{u}\tau(\mathbf{u})\leq 1-\mathbf{u}+\mathbf{u}\tau\big(\mathbf{u}^*\vee p^*\big) \leq 1-\mathbf{u}<0,
$$
since either $\tau(p^*)\leq 0$ or $\mathbf{u}^*>p^*$ and $\tau\big(\mathbf{u}^*\vee p^*\big)=0$. Thus $\mathbf{v_1}\to\infty, n\to \infty$ and, therefore, for any $v\in [1,\mathbf{v}_1]$
$$
\widetilde{\bleta}^{\beta_jr_j}_j(v,\mathbf{u})\leq \widetilde{\bleta}^{\beta_jr_j}_j(\mathbf{v_1},\mathbf{u})=\big\{\ma^{-2}\delta_n\big\}^{\frac{\mathbf{u}\omega(0)}{\mathbf{u}+\omega(0)}}
\mathbf{v}_1^{-\kappa_0(p^*,\mathbf{u})}, \quad j=1,\ldots,d.
$$
Note that
$
1-\mathbf{u}/\omega(0)+1/\beta(0)=\kappa_0(p^*,\mathbf{u})\big[1/\mathbf{u}+1/\omega(0)\big]-(\mathbf{u}-p^*)\big[1/\mathbf{u}+1/\omega(0)\big]
$
and, therefore
$$
-\frac{\kappa_0(p^*,\mathbf{u})}{1-\mathbf{u}/\omega(0)+1/\beta(0)}\geq -\frac{\mathbf{u}\omega(0)}{\mathbf{u}+\omega(0)},
$$
which yields  $\mathbf{v}_1^{-\kappa_0(p^*,\mathbf{u})}\leq \big\{\ma^{-2}\delta_n\big\}^{-\frac{\mathbf{u}\omega(0)}{\mathbf{u}+\omega(0)}}$.

It remains to note that if $\tau(\infty)\geq 0$ then $\mathbf{u}^*=\infty$ and, therefore $\mathbf{u}=\infty$. It implies  $\mathbf{v}_1=1$
and $\cI_{\mathbf{u}}(\alpha)=[\Bv,1]$ and this case has been already treated.
This completes the proof of (\ref{eq1:proof-prop2:auxiliary results-new}).

\vskip0.1cm

$\mathbf{2^0}.\;$  Remark that there obviously exists $0<\mathbf{S}:=S\big(\vec{\beta},\vec{r},\vec{\mu},p\big)<\infty$ independent of $\vec{L}$ such that
\begin{equation*}
\lim_{n\to\infty}\;(\ln{n})^{-1}\sup_{\alpha\in\{0,1\}}\sup_{\mathbf{u}\in[1,\infty]}\;
\sup_{v\in\cI_{\mathbf{u}}(\alpha)}\sum_{j=1}^d\big|\ln{\big(\blh_j(v,\mathbf{u})\big)}\big|= \mathbf{S}.
\end{equation*}
Hence, in view of (\ref{eq1:proof-prop2:auxiliary results-new}) one has for all $n$ large enough and $v\in\cI_{\mathbf{u}}(\alpha)$
\begin{eqnarray}
\label{eq4:proof-prop2:auxiliary results-new}
&&\;\;F_n\big(\vec{\blh}(v,\mathbf{u})\big)\leq\frac{\sqrt{(\mathbf{S}+2)\ln{n}}}{\sqrt{n}
\prod_{j=1}^d\big(\blh_j(v,\mathbf{u})\big)^{\frac{1}{2}+\blg_j(\alpha)}},\quad G_n\big(\vec{\blh}(v,\mathbf{u})\big)\leq\frac{(\mathbf{S}+2)\ln{n}}{n\prod_{j=1}^d\big(\blh_j(v,\mathbf{u})\big)^{1+\blg_j(\alpha)}}.
\end{eqnarray}
Taking into account that $\Bh_j(v,\mathbf{u})\geq e^{-1}\widetilde{\bleta}_j(v,\mathbf{u})$ and setting $S_0=\big[\mathbf{S}+2\big]\; e^{d+2\sum_{j=1}^d\mu_j}\BL^{-\frac{1}{\beta(1)}}$ we obtain from (\ref{eq02:proof-prop1:auxiliary results-new}) for any $\alpha\in [0,1]$ and $v\in\cI_{\mathbf{u}}(\alpha)$
\begin{gather}
\label{eq5:proof-prop2:auxiliary results-new}
(\mathbf{S}+2)n^{-1}\ln{(n)}\prod_{j=1}^d\big(\blh_j(v,\mathbf{s})\big)^{1+2\blg_j(\alpha)}\leq 2S_0\ma^2v^2\mz^{-1}(v).
\end{gather}
From now on we choose $\ma\leq a/(2S_0)<1$. It yields in view of (\ref{eq4:proof-prop2:auxiliary results-new}) and (\ref{eq5:proof-prop2:auxiliary results-new})
\begin{gather}
\label{eq7:proof-prop2:auxiliary results-new}
F_n^2\big(\vec{\blh}(v,\mathbf{u})\big)\leq a^2v^2\mz^{-1}(v),\;\;\forall v\in\cI_{\mathbf{u}}(\alpha).
\end{gather}

$\mathbf{3^0}.\;$ Since (\ref{eq7:proof-prop2:auxiliary results-new}) holds, to finish the proof of Proposition (\ref{prop2:auxiliary results-new}) all we need to show is that $G_n\big(\vec{\mh}(\mathbf{v},\mathbf{u})\big)\le a\mathbf{v},\;\;\forall v\in\cI_{\mathbf{u}}(\alpha).$ Let us distinguish three cases.

\vskip0.1cm

$\mathbf{3^0a}.\;$ Let $p^*=\infty$ or $\alpha\neq 1, \mathbf{u}=\infty$. First we note that in these cases  $\cI_{\mathbf{u}}(\alpha)=[\Bv,1]$. Next
 in view of the second inequality in (\ref{eq4:proof-prop2:auxiliary results-new}),  (\ref{eq1:proof-prop2:auxiliary results-new}), (\ref{eq5:proof-prop2:auxiliary results-new}) and (\ref{eq7:proof-prop2:auxiliary results-new}) we obtain
\begin{gather}
\label{eq70:proof-prop2:auxiliary results-new}
G_n \big(\vec{\blh}(v,\mathbf{u})\big)
\leq\frac{(\mathbf{S}+2)\ln{n}}{n\prod_{j=1}^d\big(\blh_j(v,\mathbf{u})\big)^{1+2\blg_j(\alpha)}}
\leq a^2v^2\mz^{-1}(v)\leq av,\;\;\forall v\in\cI_{\mathbf{u}}(\alpha).
\end{gather}
To get the last inequality we have used that $a<1$, $\mz(\cdot)\geq 2$ and $v\leq 1$. 

\vskip0.1cm

$\mathbf{3^0b}.\;$ Let $\alpha\neq 1, p^*<\infty, \mathbf{u}<\infty$. We have  in view of the second inequality in (\ref{eq4:proof-prop2:auxiliary results-new}) and (\ref{eq5:proof-prop2:auxiliary results-new})
$$
G_n\big(\vec{\blh}(v,\mathbf{u})\big)\leq a^2v^2\mz^{-1}(v),\quad \forall v\in\cI_{\mathbf{u}}(0).
$$
For any $\mathbf{u}\neq\infty$, simple algebra shows that $v\mz^{-1}(v)=\big\{\ma^{-2}\delta_n\big\}^{\frac{\mathbf{u}\omega(0)}{\mathbf{u}+\omega(0)}}
v^{\frac{\mathbf{u}-\omega(0)-\omega(0)/\beta(0)}{\mathbf{u}+\omega(0)}}$, and since  $\mathbf{u}\ge \mathbf{u}^*$, which implies $\mathbf{u}-\omega(0)-\omega(0)/\beta(0)>0$, the result follows from
$$
\sup_{v\in\cI_\mathbf{u}(0)}v\mz^{-1}(v)=\mathbf{v_1}\mz^{-1}(\mathbf{v_1})=1.
$$

\vskip0.1cm

$\mathbf{3^0c}.\;$ Let $\alpha=1, p^*<\infty$.
We have  in view of the second inequality in (\ref{eq4:proof-prop2:auxiliary results-new}) and (\ref{eq5:proof-prop2:auxiliary results-new})
\begin{eqnarray}
\label{eq8:proof-prop2:auxiliary results-new}
G^2_n\big(\vec{\blh}(v,\mathbf{u})\big)&\leq& \frac{(\mathbf{S}+2)\ln{n}}{n\prod_{j=1}^d\big(\blh_j(v,\mathbf{u})\big)^{1+2\blg_j(\alpha)}}
\frac{(\mathbf{S}+2)\ln{n}}{n\prod_{j=1}^d\big(\blh_j(v,\mathbf{u})\big)}
\nonumber\\
&\leq&
S_1\ma^2a^2v^2\mz^{-1}(v)\frac{\ma^{-2}\delta_n}{\prod_{j=1}^d\widetilde{\bleta}_j(v,\mathbf{u})}
,\quad v\in\cI_{\mathbf{u}}(1),
\end{eqnarray}
where we have denoted $S_1=S_0\BL^{-1/\beta(0)}$.

Our goal now is to show that for all $n$ large enough
\begin{equation}
\label{eq9:proof-prop2:auxiliary results-new}
\sup_{v\in \cI_{\mathbf{u}}(1)}\ma^{-2}\delta_n\mz^{-1}(v)\prod_{j=1}^d\widetilde{\bleta}^{-1}_j(v,\mathbf{u})\leq 1.
\end{equation}
 We easily compute  for any $v>0$
\begin{equation}
\label{eq10:proof-prop2:auxiliary results-new}
\ma^{-2}\delta_n\mz^{-1}(v)\prod_{j=1}^d\widetilde{\bleta}^{-1}_j(v,\mathbf{u})=
\mz^{-1}(v)\big\{\ma^{-2}\delta_n\big\}^{1-\frac{\mathbf{u}\omega(1)}{(\omega(1)+\mathbf{u})\omega(0)}}
\; v^{\frac{\mathbf{u}\omega(1)(2+1/\beta(1))}{(\mathbf{u}+\omega(1))\omega(0)}-\frac{1}{\beta(0)}}.
\end{equation}
Denoting the right hand side of the obtained inequality by $P(v)$ we obviously have
\begin{equation}
\label{eq11:proof-prop2:auxiliary results-new}
\sup_{v\in \cI_{\mathbf{u}}(1)}\ma^{-2}\delta_n\mz^{-1}(v)\prod_{j=1}^d\widetilde{\bleta}^{-1}_j(v,\mathbf{u})\leq
\max\big[P(\Bv),P(\mathbf{\widetilde{v}})\big],
\end{equation}
where $\mathbf{\widetilde{v}}\in\{\mathbf{v_3},\overline{\mathbf{v}},\overline{\mathbf{v}}\wedge\mathbf{v_3}\}$.
Remarking that $\mz(\Bv)=2$ we easily compute that for any $\mathbf{u}\in [1,\infty]$
\begin{equation}
\label{eq12:proof-prop2:auxiliary results-new}
P(\Bv)=2^{-1}\big\{\ma^{-2}\delta_n\big\}^{\frac{2+1/\beta(\alpha)-1/\beta(0}{2+1/\beta(\alpha)}}\to 0,\;n\to\infty.
\end{equation}
Moreover we obviously have
\begin{equation}
\label{eq13:proof-prop2:auxiliary results-new}
P(v)=2^{-1}
\big\{\ma^{-2}\delta_n\big\}^{\frac{2\mathbf{u}\omega(1)(Y+1/\mathbf{u})}{\omega(1)+\mathbf{u}}}
\; v^{\frac{2\mathbf{u}\omega(1)\pi(\mathbf{u})}{\mathbf{u}+\omega(1)}},\quad v>0.
\end{equation}

\vskip0.1cm

$\mathbf{3^0c1}.\;$ Consider the case $\kappa_1(p^*,\mathbf{u})\geq 0$. Here $\widetilde{\mathbf{v}}=\mathbf{v_3}$.

\noindent If $\pi(\mathbf{u})\leq 0$ then $\mathbf{v}_3=\infty$ and we deduce from (\ref{eq13:proof-prop2:auxiliary results-new})
\begin{equation}
\label{eq130:proof-prop2:auxiliary results-new}
\sup_{v\geq \Bv}P(v)=P(\Bv)\to 0, n\to\infty,
\end{equation}
thanks to (\ref{eq12:proof-prop2:auxiliary results-new}). If $\pi(\mathbf{u})> 0$ the definition of $\mathbf{v_3}$ implies that
\begin{equation}
\label{eq131:proof-prop2:auxiliary results-new}
P(\mathbf{v}_3)=1.
\end{equation}
Both last results together with (\ref{eq11:proof-prop2:auxiliary results-new}) and
(\ref{eq12:proof-prop2:auxiliary results-new}) prove
 (\ref{eq9:proof-prop2:auxiliary results-new})  in the case $\kappa_1(p^*,\mathbf{u})\geq 0$.

\vskip0.1cm

$\mathbf{3^0c2}.\;$ Consider the case $\kappa_1(p^*,\mathbf{u})< 0,\; Y\geq[X+1]\mathbf{y}^{-1}-1/\mathbf{u}$. Here $\widetilde{\mathbf{v}}=\overline{\mathbf{v}}$.

If $\tau(p^*)> 0$ then $\overline{\mathbf{v}}=\mathbf{v}$. Moreover
$\mathbf{y}=\mathbf{u}^*$ since $\mathbf{u}^*=\infty$ if $\tau(\infty)\geq 0$ and $\tau(\mathbf{u}^*)=0$ if $\tau(\infty)<0$.
Hence in view of
(\ref{eq900:proof-prop1:auxiliary results-new}) of Lemma \ref{lem:auxilary-results-for-prpositions-1-and2}
$$
z(1)/\omega(1)-1+2/\mathbf{u}\geq 0.
$$

We have in view of the definition of $\mathbf{v}$
\begin{eqnarray}
\label{eq14:proof-prop2:auxiliary results-new}
&&P(\mathbf{v})= 2^{-1} \big\{\ma^{-2}\delta_n\big\}^{\frac{\mathbf{u}\omega(1)(1/\omega(1)-1/\omega(0)+2/\mathbf{u})}{\omega(1)+\mathbf{u}}
+\frac{\mathbf{u}\omega^2(1)\tau(\infty)\beta(0)\pi(\mathbf{u})}{[\mathbf{u}+\omega(1)][z(\alpha)+\omega(\alpha)/\mathbf{u}]}}.
\end{eqnarray}
Note that,
\begin{eqnarray*}
&&\frac{\mathbf{u}\omega(1)(1/\omega(1)-1/\omega(0)+2/\mathbf{u})}{\omega(1)+\mathbf{u}}
+\frac{\mathbf{u}\omega^2(1)\tau(\infty)\beta(0)\pi(\mathbf{u})}{[\mathbf{u}+\omega(1)][z(1)+\omega(1)/\mathbf{u}]}
\\
&=&1-\frac{\omega(1)[1/\omega(0)-1/\mathbf{u}]}{z(1)+\omega(1)/\mathbf{u}}-\frac{\omega(1)\tau(\infty)}{z(1)+\omega(1)/\mathbf{u}}
=1-\frac{\omega(1)[1-1/\mathbf{u}]}{z(1)+\omega(1)/\mathbf{u}}>0.
\end{eqnarray*}
To get the last inequality we have used that
$$
1-\frac{\omega(1)[1-1/\mathbf{u}]}{z(1)+\omega(1)/\mathbf{u}}>0\quad\Leftrightarrow\quad z(1)/\omega(1)-1+2/\mathbf{u}>0.
$$
 Thus, we conclude that
$
P(\mathbf{v})\leq 1,
$
which together with (\ref{eq12:proof-prop2:auxiliary results-new}) implies (\ref{eq9:proof-prop2:auxiliary results-new}) in the considered case.

\vskip0.1cm

If $\tau(p^*)< 0$ then  $\overline{\mathbf{v}}=\mathbf{v_2}$. Moreover
$\mathbf{y}=\mathbf{p}^*$.  We have in view of the definition of $\mathbf{v_2}$
\begin{eqnarray}
\label{eq15:proof-prop2:auxiliary results-new}
&&P(\mathbf{v}_2)=2^{-1} \big\{\ma^{-2}\delta_n\big\}^{\frac{\mathbf{u}\omega(1)(1/\omega(1)-1/\omega(0)+2/\mathbf{u})}{\omega(1)+\mathbf{u}}
+\frac{[\mathbf{u}\omega(1)]^{2}\pi(\mathbf{u})}{\kappa_1(p^*,\mathbf{u})[\mathbf{u}+\omega(1)]^2}}.
\end{eqnarray}
After routine computations we come to the following equality
\begin{eqnarray*}
&&\frac{\mathbf{u}\omega(1)(1/\omega(1)-1/\omega(0)+2/\mathbf{u})}{\omega(1)+\mathbf{u}}
+\frac{[\mathbf{u}\omega(1)]^{2}\pi(\mathbf{u})}{\kappa_1(p^*,\mathbf{u})[\mathbf{u}+\omega(1)]^2}
\\*[2mm]
&=&-\frac{2\mathbf{u}\omega(1)p^*\big[Y-(X+1)(\mathbf{y})^{-1}+1/\mathbf{u}
\big]}
{\kappa_1(p^*,\mathbf{u})[\mathbf{u}+\omega(1)]}\geq 0.
\end{eqnarray*}
Hence, $P(\mathbf{v_2})\leq 1$ for all $n$ large enough, which together with  (\ref{eq12:proof-prop2:auxiliary results-new}) allows us to assert (\ref{eq9:proof-prop2:auxiliary results-new}) in the considered case.

\vskip0.1cm

$\mathbf{3^0c3}.\;$ Consider the case $\kappa_1(p^*,\mathbf{u})< 0,\; Y<[X+1]\mathbf{y}^{-1}-1/\mathbf{u}$. Here $\widetilde{\mathbf{v}}=\overline{\mathbf{v}}\wedge\mathbf{v_3}$.

If $\pi(\mathbf{u})\leq 0$ the required result
follows from (\ref{eq130:proof-prop2:auxiliary results-new}).

If  $\pi(\mathbf{u})> 0$ then by (\ref{eq13:proof-prop2:auxiliary results-new}) $P(\cdot)$ is strictly increasing and, therefore,
$$
P\big(\overline{\mathbf{v}}\wedge\mathbf{v_3}\big)\leq P\big(\mathbf{v_3}\big)=1
$$
in view of (\ref{eq131:proof-prop2:auxiliary results-new}). This completes the proof (\ref{eq9:proof-prop2:auxiliary results-new}).

\vskip0.1cm

Finally to conclude in the case $\mathbf{3^0c}$, choosing $\ma\leq \sqrt{1/S_1}$, we deduce from (\ref{eq8:proof-prop2:auxiliary results-new}) and (\ref{eq9:proof-prop2:auxiliary results-new})  that for all $n$ large enough

\vskip0.2cm

\centerline{$
G_n\big(\vec{\blh}(v,\mathbf{u})\big)\leq \sqrt{S_1}\ma av \leq av,
\quad v\in\cI_{\mathbf{u}}(1).
$}
%
%

\subsection{Proof of Proposition \ref{prop:measure-of-bias-deconv}}
\label{sec:subsec-Proof-prop:measure-of-bias-deconv}

In view of Lemma 5 in \cite{lepski15}, if $\tau(p^*)>0$ then
\begin{equation}
\label{eq:embedd-nik}
 \bN_{\vec{r},d}\big(\vec{\beta},\vec{L}\big) \subseteq
\bN_{\vec{q},d}\big(\vec{\gamma},c_2\vec{L}\big),
\end{equation}
where $c_2$ is independent on $\vec{L}$. Note also  that $\gamma_j\leq\beta_j$ for any $j=1,\ldots,d$.

\vskip0.1cm

$\mathbf{1^0}.\;$ Let $\big(\vec{\pi},\vec{s}\big)$ be either $\big(\vec{\beta},\vec{r}\big)$ or $\big(\vec{\gamma},\vec{q}\big)$ and without further mentioning
the couple $\big(\vec{\gamma},\vec{q}\big)$ is used below under the condition $\tau(p^*)>0$.
We obviously have for any $\vec{\mathbf{h}}\in\cH$
\begin{eqnarray*}
b_{\mathbf{h},f,j}(x)&:=&\sup_{h\in\cH:\: h\leq \mathbf{h}}
\bigg|\int_{\bR} \cK_\ell(u)\big[f\big(x+uh\mathbf{e}_j\big)-f(x)\big]\nu_1(\rd u) \bigg|
\\
&=&\sup_{h\in\cH:\: h\leq \mathbf{h}}
\bigg|\int_{\bR} \cK_\ell(u)\big[\Delta_{uh,j}f(x)\big]\nu_1(\rd u) \bigg|.
\end{eqnarray*}
For $j=1,\ldots, d$ we have
\begin{eqnarray*}
&&\int_{\bR} \cK_\ell(u) \Delta_{uh, j}f(x)\nu_1(\rd u) =
\int_{\bR} \sum_{i=1}^\ell \binom{\ell}{i} (-1)^{i+1}\frac{1}{i} \cK_\ell\Big(\frac{u}{i}\Big)\big[\Delta_{h u, j}f(x)\big]
\nu_1(\rd u)
\nonumber
\\
&&= (-1)^{\ell-1}\int_{\bR} \cK_\ell(z) \sum_{i=1}^\ell \binom{\ell}{i} (-1)^{i+\ell}\big[\Delta_{izh, j}f(x)\big]
\nu_1(\rd z)
= (-1)^{\ell-1} \int_{\bR} \cK_\ell(z) \big[\Delta^\ell_{zh, j}\, f(x)\big]\nu_1(\rd z).
\label{eq:int-representation}
\end{eqnarray*}
The last equality follows from the definition of the $\ell$-th order difference operator
(\ref{eq:Delta}).
Hence, for any $j\in\cJ_\infty$ we have in view of the definition of the Nikol'skii class (remind that $\gamma_j=\beta_j, j\in\cJ_\infty$)

$$
\|b_{\mathbf{h},f,j}\|_\infty\leq \sup_{h\in\cH:\: h\leq \mathbf{h}}\int_{\bR} \cK_\ell(z) \big\|\Delta^\ell_{zh, j}\, f(\cdot)\big\|_\infty\nu_1(\rd z)\leq L_j \sup_{h\in\cH:\: h\leq \mathbf{h}}h_j^{\pi_j} \int_{\bR} \big|\cK_\ell(z)\big||z|^{\pi_j}\nu_1(\rd z).
$$
This yields  for any $\mathbf{h}\in\cH$
\begin{equation}
\label{eq1:proof-prop3-new}
\mathbf{B}_{j,\infty,\bN_{\vec{r},d}\big(\vec{\beta},\vec{L}\big)}(\mathbf{h})\leq c_1L_j\mathbf{h}^{\pi_j},
\end{equation}
and the first  and the second assertions of the proposition are proved   for any $j\in\cJ_\infty$.

Let  $j\in\bar{\cJ}_\infty$.  Choosing
$\mathbf{k}$ from the relation $ e^{\mathbf{k}}=\mathbf{h}$ (recall that $\mathbf{h}\in\cH$), we have for any $x\in\bR^d$
\begin{eqnarray*}
b_{\mathbf{h},f,j}(x)=\sup_{k\leq \mathbf{k}}\bigg|\int_{\bR} \cK_\ell(z) \big[\Delta^\ell_{ze^{k}, j}\, f(x)\big] \nu_1(\rd z)\bigg|
=:\lim_{l\to-\infty}\;\sup_{l\leq k\leq \mathbf{k}}\bigg|\int_{\bR} \cK_\ell(z) \big[\Delta^\ell_{ze^{k}, j}\, f(x)\big] \nu_1(\rd z)\bigg|.
\end{eqnarray*}
We have  in view of monotone convergence theorem  and the triangle inequality
\begin{eqnarray*}
\big\|b_{\mathbf{h},f,j}\big\|_{s_j}&=&\lim_{l\to-\infty}\;\sup_{l\leq k\leq \mathbf{k}}\bigg\|\int_{\bR} \cK_\ell(z) \big[\Delta^\ell_{ze^{k}, j}\, f(\cdot)\big] \nu_1(\rd z)\bigg\|_{s_j}
\\
&\leq& \sum_{k=-\infty}^\mathbf{k}\bigg\|\int_{\bR} \cK_\ell(z) \big[\Delta^\ell_{ze^{k}, j}\, f(\cdot)\big] \nu_1(\rd z)\bigg\|_{s_j}.
\end{eqnarray*}
By the Minkowski inequality for integrals [see, e.g., \cite[Section 6.3]{folland}], we obtain
\begin{eqnarray*}
\label{eq4:proof-prop:measure-of-bias-deconv}
&& \big\|b_{\mathbf{v},f,j}\big\|_{s_j}
\leq\sum_{k=-\infty}^\mathbf{k}
  \int_{\bR} |\cK_\ell(z)|\big\| \Delta^\ell_{ze^{k}, j}\, f\big\|_{s_j}\;\nu_1(\rd z), \quad j=1,\ldots,d.
\end{eqnarray*}
 Taking into account that $f\in \bN_{\vec{r},d}\big(\vec{\beta},\vec{L}\big)$  and (\ref{eq:embedd-nik}), we have for any $j=1,\ldots, d,$
\begin{eqnarray}
\label{eq41:proof-prop:measure-of-bias-deconv}
\big\|b_{\mathbf{h},f,j}\big\|_{s_j}
\leq \bigg[\int_\bR |\cK_\ell(z)|\,|z|^{\beta_j} \nu_1(\rd z)\bigg]L_j\sum_{k=-\infty}^\mathbf{k} e^{k\pi_j}\leq c_1L_j\mathbf{h}^{\pi_j},\quad\forall
\mathbf{h}\in\cH^d.&&
\end{eqnarray}
This proves the first and the second assertions of the proposition for any $j\in\bar{\cJ}_\infty$.

\vskip0.1cm

$\mathbf{2^0}.\;$ Set $\bF=\bN_{\vec{r},d}\big(\vec{\beta},\vec{L}\big)$ and recall that
\vskip0.1cm
\centerline{$
\mathbf{B}^*_{j,s_j,\bF}(\mathbf{h}):=
\displaystyle{\sup_{f\in\bF}\sum_{h\in\cH:\: h\leq \mathbf{h}}}\bigg\|\int_{\bR} \cK_\ell(u)\big[f\big(x+uh\mathbf{e}_j\big)-f(x)\big]\nu_1(\rd u) \bigg\|_{s_j}\leq
\sup_{f\in\bF}\sum_{h\in\cH:\: h\leq \mathbf{h}}\big\|b_{h,f,j}\big\|_{s_j}.
$}
\vskip0.1cm
\noindent Hence, the third assertion follows from (\ref{eq1:proof-prop3-new}) and (\ref{eq41:proof-prop:measure-of-bias-deconv}).
\epr

\section{Appendix}

\subsection{Proof of Lemma \ref{lem:auxilary-results-for-prpositions-1-and2}}
Note that
\begin{eqnarray*}
z(\alpha)+\omega(\alpha)/\mathbf{s}&=&\omega(\alpha)
(2+1/\beta(\alpha))\beta(0)\tau(p^*)+1-\omega(\alpha)
(2+1/\beta(\alpha))(p^*)^{-1}+\omega(\alpha)/\mathbf{s}
\\
&=&
\omega(\alpha)(2+1/\beta(\alpha))\beta(0)\tau(p^*)-(p^*)^{-1}(1+\omega(\alpha)/\mathbf{s})\kappa_\alpha(p^*,\mathbf{s}),
\end{eqnarray*}
and (\ref{eq4:proof-prop1:auxiliary results-new}) follows. On the other hand we have
$$
z(\alpha)/\omega(\alpha)-1+2/\mathbf{u}=\big(2+2X\big)\beta(0)\tau(\infty)+2Y+2/\mathbf{u}
$$
and (\ref{eq900:proof-prop1:auxiliary results-new}) is checked if $\tau(\infty)\geq 0$ since $X,Y\geq 0$. If $\tau(\infty)<0$ and $\tau(p^*)>0$ then we note first that necessarily $\mathbf{u}^*>p^*$ since $\tau(\mathbf{u}^*)=0$ and $\tau(\cdot)$ is strictly decreasing. Hence $\mathbf{y}=\mathbf{u}^*$ and
$$
z(\alpha)/\omega(\alpha)-1+2/\mathbf{u}=\big(2+2X\big)\beta(0)\tau(\infty)+2Y+2/\mathbf{u}=2\big\{Y-(X+1)\mathbf{y}^{-1}+1/\mathbf{u}\big\}\geq 0
$$
and (\ref{eq900:proof-prop1:auxiliary results-new}) is established.

Let us prove (\ref{eq3:new-formulas}).
First we note that (\ref{eq3:new-formulas}) is obvious if $\tau(\infty)\geq 0$ because in this case
$\tau(s)>0$ for any $s\geq 1$. Thus, from now on we will assume that $\tau(\infty)<0$.

Next, if $\mathbf{u}^*> p^*$ then (\ref{eq3:new-formulas}) holds.
Indeed, in this case $0<Y-[X+1]\mathbf{y}^{-1}=Y-[X+1](\mathbf{u}^*)^{-1}$ implies $\mathbf{u}^*>(X+1)/Y$. Hence any number from the interval $\big(p^*\vee(X+1)/Y,\mathbf{u}^*\big)$ satisfies (\ref{eq3:new-formulas}).
At last,  note that if $p^*\geq \mathbf{u}^*$ we have
\begin{eqnarray*}
0\leq \frac{\kappa_\alpha(p^*,\mathbf{\infty})}{\omega(\alpha)}&=&
2+2X-2p^*Y+1/\beta(0)-p^*/\omega(0)
\\
&=&2\mathbf{y}\big[(1+X)\mathbf{y}^{-1}-Y\big]+1/\beta(0)-p^*/\omega(0)<0,
\end{eqnarray*}
since $1/\beta(0)\leq p^*/\omega(0)$ in view of $r_j\leq p^*$ for any $j=1,\ldots,d$. The obtained contradiction completes the proof of (\ref{eq3:new-formulas}).
\epr

\subsection{\textsf{Proof of Lemma \ref{lem:relations-between-parameters}}}

Indeed,
\begin{eqnarray*}
&&1/\gamma(\alpha)-1/\beta(\alpha)=1/\gamma_\pm(\alpha)-1/\beta_\pm(\alpha)=\sum_{j\in\cJ_\pm}\frac{1+2\blg_j(\alpha)}{\beta_j}\big[\tau(r_j)/\tau(p_\pm)-1\big]
\\
&&=\big[\beta(0)\tau(p_\pm)\big]^{-1}\sum_{j\in\cJ_\pm}\frac{1+2\blg_j(\alpha)}{\beta_j}(1/r_j-1/p_\pm)
=\big[\tau(p_\pm)\beta(0)\big]^{-1}
\big[1/\omega(\alpha)-1/(\beta_\pm(\alpha)p_\pm)\big].
\end{eqnarray*}
Moreover, in view of the latter inequality
\begin{eqnarray*}
1/\omega(\alpha)-1/\upsilon(\alpha)&=&1/\omega(\alpha)-1/(p_\pm\gamma_\pm(\alpha))
\\
&=&1/\omega(\alpha)-1/(p_\pm\beta_\pm(\alpha))-\big[\tau(p_\pm)\beta(0)p_\pm\big]^{-1}
\big[1/\omega(\alpha)-1/(\beta_\pm(\alpha)p_\pm)\big]
\\
&=&\big\{1-\big[\tau(p_\pm)\beta(0)p_\pm\big]^{-1}\big\}
\big[1/\omega(\alpha)-1/(\beta_\pm(\alpha)p_\pm)\big].
\end{eqnarray*}
It remains to note that $1-\big[\tau(p_\pm)\beta(0)p_\pm\big]^{-1}=\tau(\infty)/\tau(p_\pm)$ and  the lemma  follows.
\epr

\bibliographystyle{agsm}

\end{document}